%% file: articulo2.tex
\documentclass[12pt]{amsart}

\usepackage[latin1]{inputenc}
\usepackage[active]{srcltx}

\usepackage[all]{xy}

\usepackage{cases}

\usepackage{enumerate} \usepackage{float} \textwidth=6truein
\usepackage{amssymb}
\newcommand{\Z}{{\mathbb Z}} 
\newcommand{\R}{{\mathbb R}}

\newcommand{\liminv}{\operatornamewithlimits{\hbox{$\varprojlim$}}}

\renewcommand{\hom}{\operatorname{hom}\nolimits}

\newcommand{\Iso}{\operatorname{Iso}\nolimits}
\newcommand{\Hom}{\operatorname{Hom}\nolimits}

\newcommand{\Aut}{\operatorname{Aut}\nolimits}%discrete groups

\renewcommand{\ker}{\operatorname{Ker}\nolimits}
\newcommand{\im}{\operatorname{Im}\nolimits}
\newcommand{\coker}{\operatorname{Coker}\nolimits}
\newcommand{\coim}{\operatorname{Coim}\nolimits}

\renewcommand{\dim}{\operatorname{rk}\nolimits}

\newcommand{\Ab}{\operatorname{Ab}\nolimits}

\newcommand{\Ob}{\operatorname{Ob}\nolimits}

\newcommand{\A}{\ifmmode{\mathcal{A}}\else${\mathcal{A}}$\fi}
\newcommand{\B}{\ifmmode{\mathcal{B}}\else${\mathcal{B}}$\fi}
\newcommand{\C}{\ifmmode{\mathcal{C}}\else${\mathcal{C}}$\fi}
\newcommand{\D}{\ifmmode{\mathcal{D}}\else${\mathcal{D}}$\fi}
\newcommand{\G}{\ifmmode{\mathcal{G}}\else${\mathcal{G}}$\fi}
\newcommand{\I}{\ifmmode{\mathcal{I}}\else${\mathcal{I}}$\fi}
\newcommand{\J}{\ifmmode{\mathcal{J}}\else${\mathcal{J}}$\fi}
\newcommand{\K}{\ifmmode{\mathcal{K}}\else${\mathcal{K}}$\fi}
\renewcommand{\O}{\ifmmode{\mathcal{O}}\else${\mathcal{O}}$\fi}
\renewcommand{\P}{\ifmmode{\mathcal{P}}\else${\mathcal{P}}$\fi}
\newcommand{\U}{\ifmmode{\mathcal{U}}\else${\mathcal{U}}$\fi}
\newcommand{\M}{\ifmmode{\mathcal{M}}\else${\mathcal{M}}$\fi}
\newcommand{\N}{\ifmmode{\mathcal{N}}\else${\mathcal{N}}$\fi}
\newcommand{\Ss}{\ifmmode{\mathcal{S}}\else${\mathcal{S}}$\fi}
\newcommand{\T}{\ifmmode{\mathcal{T}}\else${\mathcal{T}}$\fi}
\newcommand{\Ff}{\ifmmode{\mathcal{F}}\else${\mathcal{F}}$\fi}
\newcommand{\Ll}{\ifmmode{\mathcal{L}}\else${\mathcal{L}}$\fi}
\newtheorem{Thm}{Theorem}[section]
\newtheorem{Prop}[Thm]{Proposition}

\newtheorem{Lem}[Thm]{Lemma}

\theoremstyle{definition}
\newtheorem{Defi}[Thm]{Definition}
\newtheorem{Rmk}[Thm]{Remark}
\newtheorem{Ex}[Thm]{Example}

\theoremstyle{remark}
\newtheorem*{Not}{Notation}

\theoremstyle{plain}

\newcommand{\definicio}{\stackrel{\text{def}}{=}}

\title{A method for integral cohomology of posets}

\author{Antonio D{\'\i}az Ramos}
\address{Department of Mathematical Sciences\\ King's College\\
University of Aberdeen \\ ABERDEEN AB24 3UE U.K.}
\email{a.diaz@maths.abdn.ac.uk}
\date{\today}
\begin{document}

\maketitle

\input{intro}

\input{preliminaries}

\input{local}

\input{integer}

\input{global}

\input{thing-like}

\input{morse}

\input{webb}

\input{coxeter}

%\input{quillen}

%%%%%%%%%%%%%%%%%%%%%%%%%%%%%%%%%%%%%%%%%%%%%%%%%%%%%%%%
% Bibliografia
%%%%%%%%%%%%%%%%%%%%%%%%%%%%%%%%%%%%%%%%%%%%%%%%%%%%%%%%

\input{biblio}
\end{document}

%% file: intro.tex
\section{Introduction and summary}\label{section_introduction}
Homotopy type of partially ordered sets (poset for short) play a
crucial role in algebraic topology. In fact, every space is weakly
equivalent to a simplicial complex which, of course, can be
considered as a poset. Posets also arise in more specific contexts
as homological decompositions \cite{dwyer,blo2,grodal,libman} and
subgroups complexes associated to finite groups
\cite{brown,quillen-p,bouc}. The easy structure of a poset has led
to the development of several tools to study their homotopy type,
including the remarkable Quillen's theorems \cite{quillen} and their
equivariant versions by Th\'{e}venaz and Webb \cite{G-homot}. In
spite of their apparent simplicity posets are the heart of many
celebrated problems: Webb's conjecture (proven in \cite{symonds} and
generalized in \cite{markus}), the unresolved Quillen's conjecture
on the $p$-subgroup complex (see \cite{asch-smith}), or the
fundamental Alperin's conjecture (see \cite{robinson}).

In this paper we propose a method to compute integral cohomology of
posets. This toolbox will be applicable as soon as the poset has
certain local properties. More precisely, we will require certain
structure on the category under each object of the poset. By means
of homological algebra of functors we prove that, in the presence of
these local structures, the cohomology of the poset is that of a
co-chain complex
\begin{equation}\label{equ_intro_complex}
0\rightarrow M_0\stackrel{d_0}\rightarrow
M_1\stackrel{d_1}\rightarrow M_2\stackrel{d_2}\rightarrow \ldots,
\end{equation}
where $M_n$ is free with one generator for each object of ``degree
$n$'' of the poset. If the local structure is shared at a global
level by the whole poset, further developments  show that the
cohomology of the poset is that of a co-chain complex
\begin{equation}\label{equ_intro_complex_critical}
0\rightarrow B_0\stackrel{\Omega_0}\rightarrow
B_1\stackrel{\Omega_1}\rightarrow B_2\stackrel{\Omega_2}\rightarrow
\ldots,
\end{equation}
where $B_n$ is free with one generator for each ``critical'' object
of ``degree $n$'' of the poset. We also obtain the inequalities
\begin{equation}\label{equ_intro_Morse_ineq1}
b_n\leq \dim B_n
 \end{equation}
and
\begin{equation}\label{equ_intro_Morse_ineq2}
b_n-b_{n-1}+\ldots+(-1)^n b_0\leq \dim B_n-\dim
B_{n-1}+\ldots+(-1)^n \dim B_0.
\end{equation}
where the $b_\cdot$ are the Betti numbers of the geometrical
realization of the poset.

The complex (\ref{equ_intro_complex}) applies to simplex-like
posets, i.e. posets such that the category under any object is
isomorphic to (the inclusion poset of) a simplex. The notion of
simplex-like poset is half-way between semi-simplicial complexes (as
defined originally by Eilenberg and Zilber in
\cite{semi-simplicial}) and simplicial complexes in the classical
sense. In this latter case the co-chain complex
(\ref{equ_intro_complex}) reduces to the usual simplicial co-chain
complex for integral cohomology. Notice that all the examples in the
first paragraph are simplex-like posets, but not all of them are
simplicial complexes.

Equations (\ref{equ_intro_Morse_ineq1}) and
(\ref{equ_intro_Morse_ineq2}) closely resemble weak and strong Morse
inequalities and, in fact, whenever the poset is a simplicial
complex equipped with a Morse function (see
\cite{Morse-Milnor,forman_survey}) our complex
(\ref{equ_intro_complex_critical}) is similar to the associated
Morse complex on the critical simplices. It is interesting to point
that while the Morse complex is obtained after reconstructing the
poset by homotopical gluing through critical points of the Morse
function, our complex here is directly achieved through homological
algebra.

As first application we give an alternative proof of (the
cohomological part of) Webb's conjecture for saturated fusion
systems (already proven in \cite{markus}). Further applications are
related to Coxeter groups: we prove that the Coxeter complex has the
cohomology of a sphere if the group is finite and that of a point if
the group is infinite.

The layout of the paper is as follows: Section
\ref{section_preliminaries} contains preliminaries about graded
posets and homological algebra on the category of functors. In
Section \ref{section_local} we define the local structure (``local
covering family") we require on the under categories of a poset and
study its properties. Further application of these features leads,
through Section \ref{section_integer}, to a sequence of functors to
compute the integral cohomology of a graded poset. This culminates
in the co-chain complex (\ref{equ_intro_complex}). The global
structure on the poset (``global covering family") is defined  and
used in Section \ref{section_global} to obtain co-chain complex
(\ref{equ_intro_complex_critical}) and the inequalities
(\ref{equ_intro_Morse_ineq1}) and (\ref{equ_intro_Morse_ineq2}). In
Section \ref{section_locally delta} we show how simplex-like posets
fit in this context. As an example we give a poset model of the real
projective plane. Next, Section \ref{section_morse} is devoted to
show the interplay between Morse theory and global covering
families. Finally, Webb's conjecture is proven in Section
\ref{section_webb} while Coxeter groups are treated in Section
\ref{section_coxeter}.

\textbf{Notation:}{ By the symbol $\P$ we denote a category which is
a poset (see Section \ref{section_preliminaries} below), and their
objects will be denoted by $p,q,\ldots\in \Ob(\P)$. All posets will
be graded, which means roughly that objects have an assigned degree
(see below for precise definitions). Then $\P_n=\Ob_n(\P)$ denotes
the set of objects of degree $n$. By $\Ab$ and $\Ab^\P$ we denote
the category of abelian groups and of functors from $\P$ to $\Ab$
respectively. If $S$ is a set then $|S|$ denotes the number of
elements in $S$. For a category $\C$, we denote its opposite
category by $\C^{op}$ and its realization by $|\C|$ (the geometrical
realization of the (simplicial set) nerve of the category $\C$). The
functor $c_\Z:\C\rightarrow \Ab$ is the functor which sends any
object to $\Z$ and any morphism to the identity on $\Z$.

\textbf{Acknowledgements:} I would to thank my Ph.D. supervisor A.
Viruel for his support during the development of this work, and for
reading and listening to so many different versions of it. Also,
thanks to S. Shamir and F. Xu for making comments on it.

%% file: preliminaries.tex
\section{Preliminaries}\label{section_preliminaries}
In this section we introduce some homological algebra on the abelian
category of functors $\Ab^\P$ from a given category $\P$ to the
category $\Ab$ of abelian groups. We shall work over a special kind
of categories $\P$ called graded partial ordered sets (defined
below). As well we define the cohomology of $\P$ with coefficients
$F\in \Ab^\P$ and show a ``shifting argument" to compute this
cohomology. This method will be developed systematically in the next
section in case $F$ takes free abelian groups as values.

\begin{Defi}
A \emph{poset} is a category $\P$ in which, given objects $p$ and
$p'$,
\begin{itemize}
\item there is at most one arrow $p\rightarrow p'$, and \item if
there are arrows $p\rightarrow p'$ and $p'\rightarrow p$ then
$p=p'$.
\end{itemize}
\end{Defi}
Clearly a poset defined as above is exactly the same as a set
endowed with a partial order. If $\P$ is a poset we will write
sometimes $p\leq p'$ to denote that there is an arrow $p\rightarrow
p'$, and $p<p'$ if this is the case and $p\neq p'$. To define graded
posets we first need to introduce preceding relations:
\begin{Defi}
If $\P$ is a poset and $p<p'$ then $p$ \emph{precedes} $p'$ if
$p\leq p''\leq p'$ implies that $p=p''$ or $p'=p''$.
\end{Defi}

\begin{Defi}
Let $\P$ be a poset. $\P$ is called \emph{graded} if there is a
function \linebreak $deg:\Ob(\P)\rightarrow \Z$, called the
\emph{degree function} of $\P$, which is order reversing and such
that if $p$ precedes $p'$ then $deg(p)=deg(p')+1$. If $p\in \Ob(\P)$
then $deg(p)$ is called the \emph{degree} of $p$.
\end{Defi}

Notice that for a given poset $\P$ and an object $p_0\in \Ob(\P)$
the under category $(p_0\downarrow \P)$ (see \cite{MCL}) is exactly
the full subcategory with objects $\{p|p_0\leq p\}$. We also define
\begin{Defi}\label{defi_undercategories}
Let $\P$ be a graded poset, let $p_0\in \Ob(\P)$, let $n\in \Z$ and
let $S\subseteq \Z$. Then we define $(p_0\downarrow \P)$,
$(p_0\downarrow \P)_*$, $(p_0\downarrow \P)_n$ and $(p_0\downarrow
\P)_S$ as the full subcategories of $\P$ with objects $\{p|p_0\leq
p\}$, $\{p|p_0<p\}$, $\{p|p_0\leq p, deg(p)=n\}$ and  $\{p|p_0\leq
p, deg(p)\in S\}$ respectively.
\end{Defi}

From the homotopy viewpoint, restricting to graded poset means no
loss: any topological space is weakly homotopy equivalent to a
$CW$-complex which, in turn, is homotopy equivalent to a simplicial
complex. This last can be seen as a graded poset in which the degree
function is the dimension of its simplices (more precisely and
according to our definition of order reversing degree function, the
opposite of the simplicial complex is the graded poset).

The category $\Ab^\P$, with objects the functors $F:\P\rightarrow
\Ab$ and functors the natural transformations between them, is an
abelian category in which the short exact sequences are the
object-wise ones (see \cite{maclanehom}). Because it contains enough
injectives objects (see \cite{weibel}) we can define the right
derived functors of the inverse limit functor
$\liminv:\Ab^\P\rightarrow \Ab$. For a given functor $F\in\Ab^\P$ we
define the cohomology of $\P$ with coefficients in $F$ as
$$
H^*(\P;F)={\liminv}^* F.
$$
If $F$ is the constant functor of value $M\in \Ab$ then $H^*(\P;F)$
equals the cohomology of the topological realization $|\P|$, of
$\P$, with trivial coefficients $M$. A functor $F\in\Ab^\P$ will be
called \emph{acyclic} if $H^*(\P;F)=0$ for $*>0$.

\begin{Rmk}\label{remark_conditions_on_P}
In the rest of the paper we assume the following on any graded poset
$\P$:
\begin{itemize}
\item the set $\{p|p_0\leq p\}$ is finite for any $p_0\in
\Ob(\P)$, and
\item the degree function $deg$ of $\P$ takes values
$\{\ldots,3,2,1,0\}$.
\end{itemize}
The second condition above is equivalent, by definition, to the
poset $\P$ being \emph{bounded above}, as we can always consider
translations of a degree function ($deg'=deg+c$). These conditions
will be clearly fulfilled in the applications.
\end{Rmk}

Next we introduce some acyclic objects in $\Ab^\P$.
\begin{Defi}\label{definition_F'}
Let $\P$ be a graded poset and $F\in\Ab^\P$. Then $F'$ is the
functor defined by
$$F'(p_0)=\prod_{p\in (p_0\downarrow \P)} F(p)$$
on objects $p_0\in\Ob(\P)$. For a morphism $p_1\rightarrow p_0$ the
summand $F(p)$ corresponding to $p_1\leq p$ is mapped by the
identity map to itself at the summand corresponding to $p_0\leq p$
if $p_1\leq p_0 \leq p$. Otherwise it is mapped to zero.
\end{Defi}

Notice that $F'$ is built in a similar way as enough injectives are
shown to exist in $\Ab^\P$ (see \cite[243ff.]{cohn}). The next
result summarizes some interesting properties of the functor $F'$:

\begin{Thm}\label{lemma F'properties} Let $F:\P\rightarrow \Ab$ be a
functor over a graded poset. Then the following holds:
\begin{enumerate}[(a)]
\item for each $G\in \Ab^\P$ there is a bijection
$$
\xymatrix{\Hom_{\Ab^\P}(G,F')\ar[r]^<<<<{\varphi}_<<<<{\cong} &
\prod_{p\in \Ob(\P)}
\Hom_{\Ab}(G(p),F(p)),}
$$\label{lemma_F'properties_1}
\item $\liminv F'\cong \prod_{p\in \Ob(\P)} F(p),$\label{lemma_F'properties_2}
\item $F'$ is acyclic, and\label{lemma_F'properties_3}
\item $F'$ is injective in $\Ab^\P$ if and only if $F(p)$ is injective in $\Ab$ for each $p\in\Ob(\P)$.\label{lemma_F'properties_4}
\end{enumerate}
\end{Thm}
\begin{proof}
For the first part, the bijection
$\varphi:\Hom_{\Ab^\P}(G,F')\rightarrow \prod_{p\in \Ob(\P)}
\Hom_{\Ab}(G(p),F(p))$ is given by
$$
\varphi(\mu)_p=\pi_p\eta_p,
$$
where $\pi_p:F'(p)\rightarrow F(p)$ is the projection into the
summand corresponding to $p\leq p$. The second part is consequence
of (\ref{lemma_F'properties_1}) and of the isomorphism of abelian
groups
$$
\liminv H\cong\Hom_{\Ab^\P}(\Z,H)
$$
for any $H\in\Ab^\P$, where $\Z$ is the functor $\P\rightarrow \Ab$
of constant value the integers. Part (\ref{lemma_F'properties_3}) is
proven in \cite{diaz}. There it is also proven that $F'$ is
injective if and only if $F$ takes as values injective abelian
groups.
\end{proof}

To finish this section we discuss shortly how to compute higher
limits via a ``shifting argument''. Fix $F\in \Ab^\P$ and consider
$\ker_F\in \Ab^\P$ defined by
$$
\ker_F(p_0)=\bigcap_{p\in (p_0\downarrow \P)_*} \ker
F(p_0\rightarrow p)
$$
and such that sends non-identity morphisms to zero. By Theorem
\ref{lemma F'properties}\ref{lemma_F'properties_1}) if we have a
family of maps $\{\tau_p:F(p)\rightarrow \ker_F(p)\}_{p\in \Ob(\P)}$
then there is a natural transformation $\lambda:F\Rightarrow
\ker_F'$. If $\lambda$ is object-wise injective then we obtain a
short exact sequence in $\Ab^\P$
$$
0\Rightarrow F\Rightarrow \ker_F'\Rightarrow G\Rightarrow 0,
$$
where $G$ is the object-wise co-image of $\lambda$. By Theorem
\ref{lemma F'properties}\ref{lemma_F'properties_3}) $\ker_F'$ is
acyclic, and thus the long exact sequence of the derived functors
${\liminv}^*$ gives ${\liminv}^i F={\liminv}^{i-1} G$ for $i>1$ and
$\liminv^1 F=\coim\{\liminv \ker_F'\rightarrow \liminv G\}$. The
conditions in the next definition ensure that we can build such a
natural transformation $\lambda$ which is object-wise injective:
\begin{Defi}\label{Defi_p-condensed}
Let $\P$ be a graded poset, let $F:\P\rightarrow \Ab$ be a functor
and let $n\in \Z$. We say that $F$ is \emph{$n$-condensed} if
\begin{enumerate}[(a)]
\item \label{condensed_b} $F(i)=0$ if $deg(i)<n$, and\item
\label{condensed_c} $\ker_F(i)=0$ if $deg(i)>n$.
\end{enumerate}
\end{Defi}
If the functor $F$ is $n$-condensed then we can consider the natural
transformation $\lambda:F\Rightarrow \ker'_F$ given by Theorem
\ref{lemma F'properties}\ref{lemma_F'properties_1}) for the maps
$\tau_p:F(p)\rightarrow \ker_F(p)$
\begin{numcases}{\tau_p=}
1_F(p)& if $deg(p)=n$ \nonumber \\
0     & otherwise. \nonumber
\end{numcases}
Notice that we have
\begin{numcases}{\ker_F(p)=}
F(p)& if $deg(p)=n$ \nonumber \\
0     & otherwise \nonumber
\end{numcases}
by hypotheses (\ref{condensed_b}) and (\ref{condensed_c}) of
Definition \ref{Defi_p-condensed}. Moreover the functor $\ker_F'$
takes values on objects
\begin{equation}\label{equ_expresionkerF'}
\ker_F'(p_0)=\prod_{p\in (p_0\downarrow \P)_n} F(p).
\end{equation}
The homomorphism $\lambda_p:F(p)\rightarrow \ker_F'(p)$ is given by
$$
\lambda_i=\prod_{p\in (p_0\downarrow \P)_n} F(p_0\rightarrow
p):F(p_0)\rightarrow \ker_F'(p_0)=\prod_{p\in (p_0\downarrow \P)_n}
F(p).
$$
So $\lambda_i$ is a kind of ``diagonal". An easy induction argument
on $deg(p)\in \{n,n+1,\ldots\}$ shows that $\lambda$ is a monic
natural transformation and we obtain

\begin{Lem}\label{lem_pcondensed_shortexactsequence}
Let $F:\P\rightarrow \Ab$ be an $n$-condensed functor. Then there is
a short exact sequence
$$ \xymatrix{
0\ar@{=>}[r]&F\ar@{=>}[r]^{\lambda}&\ker_F'\ar@{=>}[r]&
G\ar@{=>}[r]&0. }
$$
\end{Lem}
On the object $p_0$, $G$ takes the value
\begin{equation}\label{equ_expresionG}
G(p_0)=\prod_{p\in (p_0\downarrow \P)_n} F(p)/
\lambda_{p_0}(F(p_0)).
\end{equation}
It is clear that $G$ satisfies condition (\ref{condensed_b}) of
Definition \ref{Defi_p-condensed} for $n+1$, but in general
condition (\ref{condensed_c}) does not hold for $G$ and $n+1$. More
precisely, if $deg(p_0)>n+1$ then $\ker_G(p_0)=0$ is equivalent to
the natural map
$$
F(p_0)\rightarrow \liminv_{(p_0\downarrow \P)_*} F
$$
being an isomorphism. This natural map is a monomorphism by
condition (\ref{condensed_c}) of $F$ being $n$-condensed. So,
$\ker_G(p_0)=0$ if and only if $F(p_0)\rightarrow
\liminv_{(p_0\downarrow \P)_*} F $ is surjective. We summarize these
results in the following:

\begin{Lem}\label{lem_pcondensed_p+1condensed}
Let $F:\P\rightarrow \Ab$ be an $n$-condensed functor. Then there is
a short exact sequence
$$ \xymatrix{
0\ar@{=>}[r]&F\ar@{=>}[r]^{\lambda}&\ker_F'\ar@{=>}[r]&
G\ar@{=>}[r]&0. }
$$
Moreover, $G$ is $(n+1)$-condensed if and only if for each object
$p_0$ of degree greater than $n+1$, we have
$F(p_0)\stackrel{\cong}\rightarrow \liminv_{(p_0\downarrow \P)_*}
F$.
\end{Lem}

\begin{Ex}
Consider the graded poset $\P$ with shape
{\tiny{
$$
\xymatrix@=15pt{ \cdot_2\ar[r]\ar[rd]& \cdot_1\ar[r]\ar[rd] &\cdot_0 \\
\cdot_2\ar[r]\ar[ru]&\cdot_1\ar[r]\ar[ru]&\cdot_0,}
$$
}}where the subindexes denote the degree of the objects. The
geometrical realization $|\P|$ has the homotopy type of a two
dimensional sphere $S^2$. Consider now the functor $F:\P\rightarrow
\Ab$ with values {\tiny{
$$
\xymatrix@=15pt{ \Z_2\ar^{\times 2}[r]\ar_<<<{\times 2}[rdd]& \Z_4\ar[r]^{\times 3}\ar[rdd]^>>>{\times 6} &\Z_{12}\\
&\\
\Z_2\ar_{\times 2}[r]\ar^<<<{\times2}[ruu]&\Z_4\ar_{\times6}[r]\ar[ruu]_>>>{\times 3}&\Z_{12}.}
$$
}}
The functors $\ker_F$ and $\ker'_F$ have values
{\tiny{
$$
\xymatrix@=15pt{ 0\ar[r]^{0}\ar[rdd]_<<<{0} & 0\ar^{0}[r]\ar^>>>{0}[rdd] &\Z_{12}&&&&\Z_{12}\oplus\Z_{12}\ar[r]^{1}\ar[rdd]_<<<{1}                        &\Z_{12}\oplus\Z_{12}\ar^{1+0}[r]\ar^>>>{0+1}[rdd] &\Z_{12}\\
                         &&                                                      &&\text{and}&&    &&                  \\
 0\ar[r]_{0}\ar[ruu]^<<<{0}&0\ar_{0}[r]\ar_>>>{0}[ruu]&\Z_{12}                   &&&&\Z_{12}\oplus\Z_{12}\ar[r]_{1}\ar[ruu]^<<<{1}                        &\Z_{12}\oplus\Z_{12}\ar_{0+1}[r]\ar_>>>{1+0}[ruu] &\Z_{12}}$$
}}respectively. As $\ker_F$ is concentrated in degree $0$ the
functor $F$ is $0$-condensed. The functor $G$ from Lemma
\ref{lem_pcondensed_shortexactsequence} is given by: {\tiny{
$$
\xymatrix@=15pt{\Z_6\oplus\Z_{12}\ar^{\pi\oplus 1}[r]\ar_<<<{\pi\oplus 1}[rdd] & \Z_3\oplus\Z_{12}\ar^{0}[r]\ar^>>>{0}[rdd] &0\\
&&\\
\Z_6\oplus\Z_{12}\ar^{\pi\oplus 1}[r]\ar^<<<{\pi\oplus 1}[ruu] &\Z_3\oplus\Z_{12}\ar^{0}[r]\ar[ruu]_>>>{0}&0,}
$$}}where $\pi:\Z_6\twoheadrightarrow \Z_3$ is the quotient map. Notice
that $G$ is not $1$-condensed as $\ker_G$ takes the value $\Z_2$ on
the two objects of degree $2$ of $\P$. In fact, if $deg(p_0)=2$, then the map
$$
F(p_0)=\Z_2\rightarrow \liminv_{(p_0\downarrow \P)_*} F=\Z_4
$$
clearly is not an isomorphism. A straightforward computation shows
that
$$
\liminv F=\Z_2
$$
and
$$
{\liminv}^1 F=\coim\{\liminv \ker'_F\rightarrow \liminv
G\}=\coim\{\Z_{12}\oplus\Z_{12}\rightarrow
\Z_2\oplus\Z_6\oplus\Z_{12}\}=\Z_2.
$$
\end{Ex}

The next section is devoted to finding ``local" conditions on the
shape of the poset $\P$ such that $G$ and the subsequent functors
obtained by applying Lemma \ref{lem_pcondensed_shortexactsequence}
to a functor $F$ that takes free abelian groups as values turn out
to be condensed.

%% file: local.tex
\section{Local covering families}\label{section_local}
In this section we study a bit further the condition given in Lemma
\ref{lem_pcondensed_p+1condensed} when dealing with a functor $F$
which satisfies the following condition
\begin{Defi}\label{defi_free_functor}
Let $F:\P\rightarrow \Ab$ be a functor where $\P$ is a graded poset.
We say that $F$ is \emph{free} if $F(p)$ is a finitely generated
free abelian group for each object $p\in \Ob(\P)$ (not to be
confused with a free object in the abelian category $\Ab^\P$).
\end{Defi}

We also shall need the following

\begin{Defi}\label{defi_pure_map}
Let $A\stackrel{f}\rightarrow B$ be a map between free abelian
groups. We say that $f$ is \emph{pure} if $\coker(f)$ is a free
abelian group (see \cite{griffith}).
\end{Defi}

If $A\cong \Z^n$ is a finitely generated free abelian group we call
$\dim(A)\definicio n$. The following property of pure maps is
straightforward, and will be used repeatedly in what follows,

\begin{Lem}\label{lema_fpuremono_imply_iso}
Let $A\stackrel{f}\rightarrow B$ be a map in $\Ab$ between free
abelian groups of the same rank. If $f$ is pure and injective then
it is an isomorphism.
\end{Lem}

Now consider the condition in Lemma
\ref{lem_pcondensed_p+1condensed} again: fix $p_0$ of degree greater
than $n+1$ and consider the map given by restriction
$$
\liminv_{(p_0\downarrow \P)_*} F=\Hom_{(p_0\downarrow
\P)_*}(c_\Z,F)\rightarrow \prod_{p\in J} F(p)
$$
over a subset $J\subseteq (p_0\downarrow \P)_n$.  If this
restriction map turns out to be injective (notice that it is
injective for $J=(p_0\downarrow \P)_n$ because $F$ is $n$-condensed)
then the composition
$$F(p_0)\rightarrow \liminv_{(p_0\downarrow \P)_*} F\rightarrow
\prod_{p\in J} F(p)$$ is also injective. If $F$ is a free functor
(Definition \ref{defi_free_functor}) then both groups $F(p_0)$ and
$\prod_{p\in J} F(p)$ are free abelian groups (because we are
assuming Remark \ref{remark_conditions_on_P}). If the map
$$F(p_0)\rightarrow \prod_{p\in
J} F(p)$$ is pure then, by Lemma \ref{lema_fpuremono_imply_iso}, the
condition $\dim F(p_0)=\sum_{p\in J} \dim F(p)$ implies that this
composition is an isomorphism and so
$F(p_0)\stackrel{\cong}\rightarrow \liminv_{(p_0\downarrow \P)_*}
F$. Thus we study the subsets $J\subseteq \Ob(\P)$ that make this
restriction map a pure monomorphism:

\begin{Defi}\label{defi_covering_family}
Let $\P$ be a graded poset with degree function $deg$. A family of
subsets $\J=\{J^{p_0}_n\}_{p_0\in \Ob(\P)\text{, $0\leq n\leq
deg(p_0)$}}$ with $J^{p_0}_n\subseteq (p_0\downarrow \P)_n$ is a
\emph{local covering family} if
\begin{enumerate}[a)]
\item \label{defi_covering_a}For each $p_0$ and $0\leq n<deg(p_0)$
it holds that $\bigcup_{p\in J^{p_0}_{n+1}} (p\downarrow \P)_n =
(p_0\downarrow \P)_n$ \item \label{defi_covering_b}For each $p_0$,
$0\leq n<deg(p_0)$ and $p\in J^{p_0}_{n+1}$ it holds that
$J^p_n\subseteq J^{p_0}_n$
\end{enumerate}
\end{Defi}

Notice that the definition above does not depend on a functor
defined over the category $\P$. Also, we have
$J^{p_0}_{deg(p_0)}=\{p_0\}$ by \ref{defi_covering_a}). The next
definition states the relation we expect between a local covering
family and an $n$-condensed free functor
\begin{Defi}\label{defi_F_is_J_determined}
Let $\P$ be a graded poset, $\J$ be a local covering family and
$F:\P\rightarrow \Ab$ be an $n$-condensed free functor. We say that
$F$ is \emph{$\J$-determined} if for any object $p_0$ of degree
greater than $n+1$ the restriction map
$$
\liminv_{(p_0\downarrow \P)_*} F\rightarrow \prod_{p\in J^{p_0}_n}
F(p)
$$
is a monomorphism and the map
$$
F(p_0)\rightarrow \prod_{p\in J^{p_0}_n} F(p)
$$
is pure. If $deg(p_0)=n+1$ then we require the last map above to be
a pure monomorphism.
\end{Defi}

The main feature of local covering families is that they allow
freeness plus $\J$-determinacy to pass from $F$ to $G$. For an
object $p_0$ with $deg(p_0)\geq n+1$ notice that the map
$$
F(p_0)\rightarrow \prod_{p\in J^{p_0}_n} F(p)
$$
is a pure monomorphism as consequence of Definition
\ref{defi_F_is_J_determined}. The condition of Definition
\ref{defi_F_is_J_determined} for $deg(p_0)=n+1$ is added in order to
obtain that $G$ is a free functor. Notice that the following
proposition restricts to functors which take free abelian groups as
values.

\begin{Prop}\label{prop_covering_familiy_determination_hereditary}
Let $\P$ be a graded poset and $\J$ a local covering family. Assume
that $F:\P\rightarrow \Ab$ is $n$-condensed, free and
$\J$-determined  and consider the functor $G$ defined by
$$ \xymatrix{
0\ar@{=>}[r]&F\ar@{=>}[r]^{\lambda}&\ker_F'\ar@{=>}[r]&
G\ar@{=>}[r]&0. }
$$
If for each object $p_0$ with $deg(p_0)\geq n+1$ it holds that $\dim
F(p_0)=\sum_{p\in J^{p_0}_n} \dim F(p)$, then $G$ is
$(n+1)$-condensed, free and $\J$-determined.
\end{Prop}
\begin{proof}
Notice that the hypothesis implies that for any object $p_0$ of
degree $deg(p_0)>n+1$ the two maps
$$
F(p_0)\rightarrow \liminv_{(p_0\downarrow \P)_*} F\rightarrow
\prod_{p\in J^{p_0}_n} F(p)
$$
are isomorphisms. In particular, $F(p_0)\stackrel{\cong}\rightarrow
\liminv_{(p_0\downarrow \P)_*} F$ and so $G$ is $(n+1)$-condensed.

If $deg(p_0)=n+1$ then the map
$$
F(p_0)\rightarrow \prod_{p\in J^{p_0}_n} F(p)
$$
is an isomorphism by hypothesis. Next we prove that $G$ is a free
functor. Consider any $p\in \Ob(\P)$ with $deg(p)\geq n+1$ (if
$deg(p)<n+1$ then $G(p)=0$) and the short exact sequence of abelian
groups
$$
0\rightarrow F(p)\stackrel{\lambda_p}\rightarrow
\ker_F'(p)\stackrel{\pi_p}\rightarrow G(p)\rightarrow 0.
$$
Then it is straightforward that the map
$$
s_p:\ker_F'(p)=\prod_{q\in (p\downarrow \P)_n}
F(q)\twoheadrightarrow \prod_{q\in J^p_n}
F(q)\stackrel{\cong}\rightarrow F(p)
$$
is a section of $\lambda_p$, i.e. $s_p\circ \lambda_p=1_{F(p)}$
(since the restriction map $F(p)\rightarrow \prod_{q\in J^{p}_n}
F(q)$ is injective). This implies that the short exact sequence
above splits and so $G(p)$ is a subgroup of the free abelian group
$\ker_F'(p)$, and thus it is free as well.

Next we prove that $G$ is $\J$-determined. Take $p_0$ of degree
$n=deg(p_0)$ greater than $n+2$. We first check that the restriction
map
$$
\liminv_{(p_0\downarrow \P)_*} G\rightarrow \prod_{p\in
J^{p_0}_{n+1}} G(p)
$$
is injective. Consider any element $\psi\in \liminv_{(p_0\downarrow
\P)_*} G=\Hom_{\Ab^\P}(\Z,G)$ which is in the kernel of the
restriction map above. Notice that, as $deg(p_0)>n+2$, we can
consider the subset $J^{p_0}_{n+2}\subseteq (p_0\downarrow \P)_*$.
If for any $q\in J^{p_0}_{n+2}$ it holds that $\psi_q(1)=0$ then
$\psi=0$ because of Definition
\ref{defi_covering_family}\ref{defi_covering_a}) and because $G$ is
$(n+1)$-condensed.

Thus take $q\in J^{p_0}_{n+2}$. We want to see that $x\definicio
\psi_q(1)=0$. Recall the short exact sequence of abelian groups
$$
0\rightarrow F(q)\stackrel{\lambda_q}\rightarrow
\ker_F'(q)\stackrel{\pi_q}\rightarrow G(q)\rightarrow 0
$$
and take $y\in \ker_F'(q)$ such that $\pi_q(y)=x$. Recall that
$\ker_F'(q)=\prod_{p\in (q\downarrow \P)_{n}} F(p)$ and denote by
$\alpha_p:q\rightarrow p$ the unique arrow from $q$ to $p$ for $p\in
(q\downarrow \P)_n$.

Now consider the restriction $y|_{J^q_n}\in \prod_{p\in J^q_n}
F(p)$. Because $deg(q)=n+2>n+1$ the map $F(q)\rightarrow \prod_{p\in
J^q_n} F(p)$ is an isomorphism by hypothesis. Then there exists a
unique $z\in F(q)$ with $F(\alpha_p)(z)=y_p$ for each $p\in
J^q_n\subseteq (q\downarrow \P)_{n}$. If we prove that
$F(\alpha_p)(z)=y_p$ for each $p\in (q\downarrow \P)_{n}$ then
$\lambda_q(z)=y$. This implies that
$x=\pi_q(y)=\pi_q(\lambda_q(z))=0$ and completes the proof.

Thus take $p\in (q\downarrow \P)_n$. By Definition
\ref{defi_covering_family}\ref{defi_covering_a}) there is
$\beta_p:p'\rightarrow p$ with $p'\in J^q_{n+1}$. Write
$\beta_{p'}:q\rightarrow p'$ for the unique arrow from $q$ to $p'$.
It holds that $\alpha_p=\beta_p\circ\beta_{p'}$. By Definition
\ref{defi_covering_family}\ref{defi_covering_b}) we have that
$J^q_{n+1}\subseteq J^{p_0}_{n+1}$. Thus
$G(\beta_{p'})(x)=G(\beta_{p'})(\psi_q(1))=\psi_{p'}(1)=0$ as $\psi$
is in the kernel of the restriction map. The short exact sequence
$$ 0\rightarrow F(p')\stackrel{\lambda_{p'}}\rightarrow
\ker_F'(p')\stackrel{\pi_{p'}}\rightarrow G(p')\rightarrow 0
$$
implies that there exists $t_{p'}\in F(p')$ such that
$\lambda_{p'}(t_{p'})=\ker_F'(\beta_{p'})(y)$. Consider
$z_{p'}=F(\beta_{p'})(z)$. We have that $z_{p'}$ and $t_{p'}$ have
the same image by the restriction map
$$
\liminv_{\P^{p'}_*} F\rightarrow \prod_{p\in J^{p'}_{n}} F(p)
$$
because $J^{p'}_n\subseteq J^q_n$. Because $F$ is $\J$-determined
then this restriction map is a monomorphism and so $z_{p'}=t_{p'}$.
This implies that
$$F(\alpha_p)(z)=F(\beta_p\circ\beta_{p'})(z)=F(\beta_p)(z_{p'})=F(\beta_p)(t_{p'})=y_p$$
and the proof of the restriction map being injective is finished.

Now we check that the map
$$
\omega:G(p_0)\rightarrow \prod_{p\in J^{p_0}_{n+1}} G(p)
$$
is pure. Take $z\in \prod_{p\in J^{p_0}_{n+1}} G(p)$ such that there
exists $x\in G(p_0)$ with $m\cdot z=\omega(x)$ for some $m\neq 0$.
We have to check that there exists $x'\in G(p_0)$ with
$z=\omega(x')$, or equivalently, that $x=m\cdot x'$ for some $x'\in
G(p_0)$. Recall once more the short exact sequence of abelian groups
$$
0\rightarrow F(p_0)\stackrel{\lambda_{p_0}}\rightarrow
\ker_F'(p_0)\stackrel{\pi_{p_0}}\rightarrow G(p_0)\rightarrow 0
$$
and take $y\in \ker_F'(p_0)$ with $\pi_{p_0}(y)=x$. We are going to
build $h\in F(p_0)$ such that $y-\lambda_{p_0}(h)=m\cdot y'$, i.e.,
such that for any $p\in (p_0\downarrow \P)_n$ the element
$(y-\lambda_{p_0}(h))_p=y_p-F(p_0\rightarrow p)(h)\in F(p)$ is
divisible by $m$. This implies that $x=m\cdot x'$ with
$x'=\pi_{p_0}(y')$.

Notice that by hypothesis for each $q\in J^{p_0}_{n+1}$,
$G(p_0\rightarrow q)(x)=m\cdot z_q \in G(q)$. This implies that
there exist $h_q\in F(q)$ and $y_q\in \ker_F'(q)$ such that
$\ker_F'(p_0\rightarrow q)(y)-\lambda_q(h_q)=m\cdot y_q$, i.e., such
that for each $p\in (q\downarrow \P)_n\subseteq (p_0\downarrow
\P)_n$  we have that $y_p-F(q\rightarrow p)(h_q)=m\cdot (y_q)_p\in
F(p)$ (it is enough to take $y_q$ with $\pi_q(y_q)=z_q$).

To build $h$ we use the map
$$
\tau:\prod_{p\in J^{p_0}_n} F(p)\stackrel{\cong}\rightarrow F(p_0)
$$
given by hypothesis, which is the inverse of the map
$$
F(p_0)\rightarrow \prod_{p\in J^{p_0}_n} F(p).
$$
For each $p\in J^{p_0}_n\subseteq (p_0\downarrow \P)_n$ choose, by
Definition \ref{defi_covering_family}\ref{defi_covering_a}),
$q(p)\in J^{p_0}_{n+1}$ such that there is an arrow $q(p)\rightarrow
p$. Then set $\eta_p=F(q(p)\rightarrow p)(h_{q(p)})\in F(p)$, where
$h_{q(p)}$ is built as before. Define $h\definicio \tau(\eta)$. By
construction $F(p_0\rightarrow p)(h)=F(q(p)\rightarrow p)(h_{q(p)})$
for each $p\in J^{p_0}_n$ (but not for an arbitrary $p\in
(p_0\downarrow \P)_n$).

With this definition for $h$ we check now that $y_p-F(p_0\rightarrow
p)(h)$ is divisible by $m$ for each $p\in (p_0\downarrow \P)_n$.
This finishes the proof. Fix $p\in (p_0\downarrow \P)_n$ and $q_p\in
J^{p_0}_{n+1}$ such that there is an arrow $q_p\rightarrow p$ (we
are not assuming that $q_p=q(p)$ if $p\in J^{p_0}_n$). On the one
hand we have by hypothesis that
$$
y_k-F(q_p\rightarrow k)(h_{q_p})=m\cdot (y_{q_p})_k
$$
for each $k\in \P^{q_p}_n$. In particular,
\begin{equation}\label{equ_covefami_1}
y_k-F(q_p\rightarrow k)(h_{q_p})=m\cdot (y_{q_p})_k
\end{equation}
for each $k\in J^{q_p}_n$. Set $h'\definicio F(p_0\rightarrow
q_p)(h)$. Because $q_p\in J^{p_0}_{n+1}$ then, by Definition
\ref{defi_covering_family}\ref{defi_covering_b}),
$J^{q_p}_n\subseteq J^{p_0}_n$ and thus by construction for any
$k\in J^{q_p}_n$
$$
y_k-F(p_0\rightarrow k)(h)=y_k-F(q(k)\rightarrow k)(h_{q(k)})=m\cdot
(y_{q(k)})_k.
$$
Notice that $F(p_0\rightarrow k)(h)=F(q_p\rightarrow
k)(F(p_0\rightarrow q_p)(h))=F(q_p\rightarrow k)(h')$. On the other
hand, we have obtained
\begin{equation}\label{equ_covefami_2}
y_k-F(q_p\rightarrow k)(h')=m\cdot (y_{q(k)})_k
\end{equation}
for each $k\in \J^{q_p}_n$.

Now write $\eta_k=(y_{q(k)})_k-(y_{q_p})_k$ for each $k\in
J^{q_p}_n$ and write $h''=\tau(\eta)\in F(q_p)$ where $\tau$ is the
inverse of the map
$$
F(q_p)\rightarrow \prod_{p\in J^{q_p}_n} F(p).
$$
By Equations (\ref{equ_covefami_1}) and (\ref{equ_covefami_2}) it is
straightforward that the elements $h_{q_p}-m\cdot h''$ and $h'$ have
the same image by this map. Then, as this map is injective by
hypothesis, $h'=h_{q_p}-m\cdot h''$. As $p\in J^{q_p}_n$ we have
$$
y_p-F(p_0\rightarrow p)(h)=y_p-F(q_p\rightarrow
p)(h')=y_p-F(q_p\rightarrow p)(h_{q_p}-m\cdot h''),
$$
and this equals
$$
m\cdot (y_{q_p})_p+m\cdot F(q_p\rightarrow p)(h'').
$$
Thus $y_p-F(p_0\rightarrow p)(h)$ is divisible by $m$.

If $deg(p_0)=n+2$ we have to see that the map
$$
\omega:G(p_0)\rightarrow \prod_{p\in J^{p_0}_{n+1}} G(p)
$$
is a pure monomorphism. To prove that $\omega$ is a monomorphism use
the proof above starting where $\psi_q$ is considered for an
arbitrary object $q$ of degree $n+2$. The proof of $\omega$ being
pure is exactly the same as above.
\end{proof}

\begin{Rmk}\label{rmk_dimensions}
Notice that in the conditions of the proposition we have the
following formula for the rank of the free abelian group $G(p_0)$
for $deg(p_0)\geq n+1$
$$
\dim(G(p_0))=\sum_{p\in (p_0\downarrow \P)_{n}} \dim F(p) - \dim
F(p_0).
$$
This is so because of the short exact sequence of free abelian
groups
$$
0\rightarrow F(p_0)\stackrel{\lambda_{p_0}}\rightarrow
\ker_F'(p_0)\stackrel{\pi_{p_0}}\rightarrow G(p_0)\rightarrow 0.
$$
\end{Rmk}
\begin{Rmk}\label{rmk_section} In the conditions of the proposition
there are isomorphisms
$$
F(p_0)\stackrel{\cong}\rightarrow \prod_{p\in J^{p_0}_n} F(p)
$$
for each $p_0$ with $deg(p_0)\geq n+1$. Moreover, we have built a
map $s_{p_0}:\ker'_F(p_0)\rightarrow F(p_0)$ with $s_{p_0}\circ
\lambda_{p_0}=1_{F(p_0)}$. To the homomorphism $s_{p_0}$ corresponds
the monomorphism
$$
\xymatrix{ \ar[r]^{\delta_{p_0}}G(p_0) & \ker'_F(p_0)}
$$
given by
$$
\xymatrix{ \pi_{p_0}(x)\ar@{|->}[r] & x-(\lambda_{p_0}\circ
s_{p_0})(x)},
$$
which satisfies $\pi_{p_0}\circ\delta_{p_0}=1_{G(p_0)}$. It is
straightforward that, by construction, $\im \delta_{p_0}=\prod_{p\in
(p_0\downarrow \P)_n\setminus J^{p_0}_n} F(p)$, and thus
$$
G(p_0)\stackrel{\delta_{p_0}}\cong \prod_{p\in (p_0\downarrow
\P)_n\setminus J^{p_0}_n} F(p).
$$
Moreover, $x=\delta_{p_0}(y)$ is the only preimage of $y$ by
$\pi_{p_0}$ which verifies $x_p=0$ for $p\in J^{p_0}_n$.
\end{Rmk}

The main consequence of the previous proposition is that it reduces
the problem of whether $G$ is $(n+1)$-condensed to some integral
equations. Moreover, this procedure can be applied recursively
because the resulting functor $G$ turns out to be $(n+1)$-condensed,
free and $\J$-determined, and so the proposition applies to $G$ too.
Notice again that the ranks of $G$ are given by Remark
\ref{rmk_dimensions}.

\begin{Ex}
Consider the graded poset $\P$ with shape {\tiny{
$$
\xymatrix@=15pt{a_2 \ar[r]\ar[rd]\ar[rdd] & c_1 \ar[r]\ar[rd]\ar[rdd]& f_0\\
&d_1\ar[r]\ar[rd]\ar[ru]&g_0\\
b_2 \ar[r]\ar[ru]\ar[ruu] & e_1 \ar[r]\ar[ru]\ar[ruu]& h_0,}
$$}}where the subindexes denote the degree of the objects. Notice that
$\P$ has the homotopy type of a wedge of $4$ $2$-dimensional spheres
$\bigvee_{i=1}^{4} S^2_i$. The sets
\begin{eqnarray}
&J^f_0=\{f\},J^g_0=\{g\},J^h_0=\{h\},\nonumber\\
&J^c_1=\{c\},J^c_0=\{f\},J^d_1=\{d\},J^d_0=\{g\},J^e_1=\{e\},J^e_0=\{h\},\nonumber\\
&J^a_2=\{a\},J^a_1=\{c\},J^a_0=\{f\},J^b_2=\{b\},J^b_1=\{e\},J^b_0=\{h\}.\nonumber
\end{eqnarray}
define a local covering family $\J$ for $\P$. Consider now the the
functor $F:\P\rightarrow \Ab$ with values {\tiny{
$$
\xymatrix@=15pt{\Z \ar[r]\ar[rd]\ar[rdd] & \Z \ar[r]\ar[rd]\ar[rdd]& \Z\\
&\Z\ar[r]\ar[rd]\ar[ru]&\Z\\
\Z \ar[r]\ar[ru]\ar[ruu] & \Z \ar[r]\ar[ru]\ar[ruu]& \Z,}
$$}}such that all the arrows arriving and departing from $c$ and $e$ are
the identity and all the arrows arriving and departing from $d$ are
minus the identity. Then $F$ is $0$-condensed, free and
$\J$-determined. Moreover, equations in Proposition
\ref{prop_covering_familiy_determination_hereditary} hold for any
object of degree greater or equal to $1$. Thus we obtain a functor
$G$ which is $1$-condensed, free and $\J$-determined. By Remark
\ref{rmk_dimensions} we know the ranks of the free abelian groups
that $G$ takes as values: {\tiny{
$$
\xymatrix@=15pt{\Z^2 \ar[r]\ar[rd]\ar[rdd] & \Z^2 \ar[r]\ar[rd]\ar[rdd]& 0\\
&\Z^2\ar[r]\ar[rd]\ar[ru]&0\\
\Z^2 \ar[r]\ar[ru]\ar[ruu] & \Z^2 \ar[r]\ar[ru]\ar[ruu]& 0.}
$$}}

Moreover, equations in Proposition
\ref{prop_covering_familiy_determination_hereditary} applied to $G$
hold for any object of degree greater or equal to $2$, and thus we
obtain a functor $H$ which is $2$-condensed, free and
$\J$-determined: {\tiny{
$$
\xymatrix@=15pt{\Z^4 \ar[r]\ar[rd]\ar[rdd] & 0 \ar[r]\ar[rd]\ar[rdd]& 0\\
&0\ar[r]\ar[rd]\ar[ru]&0\\
\Z^4 \ar[r]\ar[ru]\ar[ruu] & 0 \ar[r]\ar[ru]\ar[ruu]& 0,}
$$}}

By Remark \ref{rmk_section} we can identify $G(a)\cong F(g)\oplus
F(h)=\Z_g\oplus \Z_h$. Also $G(c)\cong \Z_g\oplus \Z_h$, $G(d)\cong
\Z_f\oplus \Z_h$ and $G(e)\cong \Z_f\oplus \Z_g$. By the definition
of $G$ as a co-image it is easy to see that
\begin{eqnarray}
G(a\rightarrow c):&\Z_g\oplus \Z_h\rightarrow\Z_g\oplus \Z_h\nonumber \\
&(g,h)\mapsto (g,h).\nonumber
\end{eqnarray}
Also $G(a\rightarrow d)(g,h)=(-g,h-g)$ and $G(a\rightarrow
e)(g,h)=(-h,g-h)$ with respect to the ordered basis mentioned above.
Additional computations lead to
$$
\liminv F=\Z,
$$
$$
{\liminv}^1 F=\coim\{\liminv \ker'_F\rightarrow \liminv
G\}=\coim\{\Z^3\rightarrow\Z^2\}=0\text{, and}
$$
$$
{\liminv}^2 F=\coim\{\liminv \ker'_G\rightarrow \liminv
H\}=\coim\{\Z^6\rightarrow\Z^8\}=\Z^4.
$$
\end{Ex}

%section "Local work: sequences of functors"
%\input{coho_local1}

%section "Local work: generators and relations"
%\input{coho_local2}

%% file: integer.tex
\section{Integral cohomology}\label{section_integer}
In this section we apply the work developed in the preceding
sections to compute the cohomology with integer coefficients of the
realization of a graded poset $\P$ equipped with a local covering
family $\J$.

To compute the abelian group $H^n(|\P|;\Z)$ for $n\geq 0$ we compute
the higher limit ${\liminv}^n c_\Z$ where $c_\Z:\P\rightarrow \Ab$
is the functor of constant value $\Z$ which sends every morphism to
the identity $1_\Z$. We begin studying the extra conditions that the
local covering family $\J$ must satisfy to apply iteratively the
Proposition \ref{prop_covering_familiy_determination_hereditary}
beginning on $c_\Z$.

First, notice that $c_\Z$ is $0$-condensed (we are assuming
\ref{remark_conditions_on_P}) and free (Definition
\ref{defi_free_functor}). By Definition
\ref{defi_F_is_J_determined}, $c_\Z$ is $\J$-determined as
$0$-condensed functor if and only if for each $p_0\in \Ob(\P)$ with
$deg(p_0)\geq 2$ the set $J^{p_0}_0$ intersects each connected
component of $(p_0\downarrow \P)_*$. The dimensional equation in
Proposition \ref{prop_covering_familiy_determination_hereditary} for
$p_0\in \Ob(\P)$ with $deg(p_0)\geq 1$ becomes $\dim
c_\Z(p_0)=1=|J^{p_0}_0|=\sum_{p\in J^{p_0}_0} \dim c_\Z(p)$. Thus,
$c_\Z$ is $\J$-determined as a $0$-condensed functor if and only if
$(p_0\downarrow \P)_*$ is connected for $deg(p_0)\geq 2$ and
$|J^{p_0}_0|=1$ for $deg(p_0)\geq 1$. The successive applications of
Proposition \ref{prop_covering_familiy_determination_hereditary}
give, by the dimensional equation in the statement of the
Proposition \ref{prop_covering_familiy_determination_hereditary},
the following:

\begin{Defi}\label{defi_structural numbers of graded poset}
Let $\P$ be a graded poset. Define, inductively on $n$, the number
$R^{p_0}_n$ for each object $p_0$ with $deg(p_0)\geq n$ by
$R^{p_0}_0=1$ and by
$$R^{p_0}_n=\sum_{p\in (p_0\downarrow \P)_{n-1}} R^p_{n-1}-R^{p_0}_{n-1}$$
for $n\geq 1$.
\end{Defi}
\begin{Defi}\label{defi_adequate_covering_family}
Let $\P$ be a graded poset and $\J$ be a local covering family for
$\P$. We say that $\J$ is \emph{adequate} if $(p_0\downarrow \P)_*$
is connected for $deg(p_0)\geq 2$, and if we have the equality
$$
R^{p_0}_n=\sum_{p\in J^{p_0}_n} R^p_n
$$
for $n\geq 0$ and $deg(p_0)\geq n+1$.
\end{Defi}

\begin{Prop}\label{cohomology with adequate local covering family}
Let $\P$ be a graded poset and let $\J$ be an adequate local
covering family. Then there is a sequence of functors
$F_0,F_1,F_2,\ldots$ defined by $F_0\definicio c_\Z:\P\rightarrow
\Ab$ and by the short exact sequence
$$ \xymatrix{
0\ar@{=>}[r]&F_{n-1}\ar@{=>}[r]^{\lambda_{n-1}}&\ker_{F_{n-1}}'\ar@{=>}[r]^{\pi_n}&
F_n\ar@{=>}[r]&0 }
$$
for $n=1,2,3,\ldots$. Moreover, $F_n$ is $n$-condensed, free and
$\J$-determined for $n\geq 0$. For $deg(p_0)\geq n$ we have $\dim
F_n(p_0)=R^{p_0}_n$.
\end{Prop}

The local properties of a graded poset $\P$ equipped with an
adequate local covering family $\J$ give rise to a sequence
$F_0=c_\Z$, $F_1$, $F_2$,\ldots of functors. Now, we study some
properties of these functors which are independent of the local
covering family $\J$.

The first point to notice is that the sequence of functors from
Proposition \ref{cohomology with adequate local covering family}
does not depend on the adequate local covering family $\J$. Thus,
two or more adequate local covering families can be considered for
the the same graded poset and they still give rise to the same
sequence of functors. Next we focus on the short exact sequence
$$ \xymatrix{
0\ar@{=>}[r]&F_{n-1}\ar@{=>}[r]^{\lambda_{n-1}}&\ker_{F_{n-1}}'\ar@{=>}[r]^{\pi_n}&
F_n\ar@{=>}[r]&0 }
$$
for $n\geq 1$ of Proposition \ref{cohomology with adequate local
covering family}.  The beginning of the long exact sequence of
derived functors for this short exact sequence is
\begin{equation}\label{eqn_long_exact_sequence_H^n}
\xymatrix{0\ar[r]&\liminv F_{n-1}\ar[r]^{\iota_{n-1}}&\liminv
\ker_{F_{n-1}}'\ar[r]^{\omega_n}& \liminv F_n\ar[r]&H^n(\P;\Z)
\ar[r]&0 }
\end{equation}
by Theorem \ref{lemma F'properties}, where
$\iota_{n-1}=\widetilde{\lambda_{n-1}}$ and
$\omega_n=\widetilde{\pi_n}$ are the induced maps. Notice that the
three inverse limits appearing above are free abelian groups as the
corresponding functors take free abelian groups as values. In fact,
for the middle term we have the exact description
\begin{equation}\label{eqn_kerF' directproduct}
\liminv \ker_{F_{n-1}}' \cong\prod_{p\in \Ob_{n-1}(\P)} F_{n-1}(p)
\end{equation}
given by Theorem \ref{lemma F'properties}. It turns out that there
is also a simpler description for $\liminv F_n$, which can be
interpreted as the analogue in the context of $CW$-complexes to the
fact that the cohomology on degree $n$ depends upon the $n+1$
skeleton (recall that $F_n(p)=0$ if $deg(p)<n$):

\begin{Lem}\label{lemma_liminvFp n+2_esqueleto}
Let $\P$ be a graded poset and let $\J$ be an adequate local
covering family. Let $c_\Z,F_1,F_2,\ldots$ be the sequence of
functors given by Proposition $\ref{cohomology with adequate local
covering family}$. Then
$$
\liminv F_n\cong\liminv F_n|_{\P_{\{n+1,n\}}}
$$
for each $n\geq 0$.
\end{Lem}
\begin{proof}
Consider the restriction map
$$
\liminv F_n\rightarrow \liminv F_n|_{\P_{\{n+1,n\}}}.
$$
This map is clearly a monomorphism because $F_n$ is an $n$-condensed
functor. To see that it is surjective take $\psi\in \liminv
F_n|_{\P_{\{n+1,n\}}}$. We want to extend $\psi$ to each $p\in
\Ob(\P)$ with $deg(p)>n+1$. We do it inductively on $deg(p)$.

Notice that (see the beginning of the proof of Proposition
\ref{prop_covering_familiy_determination_hereditary})
$$
F_n(p)\rightarrow \liminv_{(p\downarrow \P)_*} F_n
$$
is an isomorphism for $deg(p)>n+1$. For $deg(p)=n+2$ we have that
$j\in {(p\downarrow \P)_*}$ implies that $deg(j)\leq n+1$. Then
there is a unique way of extending $\psi$ to $\psi(p)$. Once we have
extended $\psi$ to $\P_{\{n+2,n+1,n\}}$ we proceed with an induction
argument. That the extension that we are building is actually a
functor is again due to that $F_n$ is $n$-condensed.
\end{proof}

Also, from Equations (\ref{eqn_long_exact_sequence_H^n}) and
(\ref{eqn_kerF' directproduct}), we have the following formula,
analogue to that of the Euler characteristic for $CW$-complexes:
\begin{Lem}\label{lemma_analogue eulerchar}
Let $\P$ be a poset for which exists an adequate local covering
family. Then
$$
\sum_i (-1)^i\dim H^i(\P;\Z)=\sum_i (-1)^i\sum_{p\in
\Ob_i(\P)}R^p_i.
$$
\end{Lem}

Take up again Equations (\ref{eqn_long_exact_sequence_H^n}) and
(\ref{eqn_kerF' directproduct}). We can form the sequence of abelian
groups
$$
0\rightarrow \prod_{p\in \Ob_0(\P)} F_0(p)\stackrel{d_0}\rightarrow
\prod_{p\in \Ob_1(\P)} F_1(p) \stackrel{d_1}\rightarrow \prod_{p\in
\Ob_2(\P)} F_2(p)\stackrel{d_2}\rightarrow \ldots
$$
where $d_n=\iota_{n+1}\circ w_{n+1}$ for $n\geq 0$. Then it is
straightforward that this sequence is a co-chain complex and its
cohomology is exactly the cohomology of $\P$ with integers
coefficients:
\begin{Thm}\label{proposition_cohomology_cochain_covering}
Let $\P$ be a graded poset for which exists an adequate local
covering family. Then there exists a co-chain complex
$$
0\rightarrow \prod_{p\in \Ob_0(\P)} F_0(p)\stackrel{d_0}\rightarrow
\prod_{p\in \Ob_1(\P)} F_1(p) \stackrel{d_1}\rightarrow \prod_{p\in
\Ob_2(\P)} F_2(p)\stackrel{d_2}\rightarrow \ldots
$$
of which cohomology is $H^*(\P;\Z)$.
\end{Thm}

%% file: global.tex
\section{Global covering families}\label{section_global}
Recall that local covering families were defined as subsets of the
local categories $(p\downarrow \P)$ for $p\in \Ob(\P)$, where $\P$
is a graded poset. In this section we define global covering
families by subsets of the whole category $\P$, imitating some of
the local features of the local covering families.

\begin{Defi}\label{defi_global covering family}
Let $\P$ be a graded poset for which there exists an adequate local
covering family, and consider the sequence of functors $F_0=c_\Z$,
$F_1$, $F_2$,\ldots given by Proposition \ref{cohomology with
adequate local covering family}. A \emph{global covering family} is
a family of subsets $\K=\{K_n\}_{n\geq 0}$ with $K_n\subseteq
\Ob_n(\P)$ such that
\begin{enumerate}
\item\label{def_globalcovfam_mono}the morphism
$$
\liminv F_n\rightarrow \prod_{p\in K_n} F_n(p)
$$
is injective for each $n\geq 0$, and
\item\label{def_globalcovfam_pure}the morphism
$$
\prod_{p\in \Ob_{n-1}(\P)\setminus K_{n-1}} F_{n-1}(p)\rightarrow
\prod_{p\in K_n} F_n(p)
$$
is pure for each $n\geq 1$.
\end{enumerate}
\end{Defi}

Notice that the definition does not depend on the particular
adequate local covering family used to obtain the sequence of
functors (as the sequence of functors does not depend on it). The
map
$$
\prod_{p\in \Ob_{n-1}(\P)\setminus K_{n-1}} F_{n-1}(p)\rightarrow
\prod_{p\in K_n} F(p)
$$
used in the definition is the composition
$$
\prod_{p\in \Ob_{n-1}(\P)\setminus K_{n-1}} F_{n-1}(p)
\hookrightarrow\prod_{p\in \Ob_{n-1}(\P)}
F_{n-1}(p)\stackrel{d_{n-1}}\rightarrow \prod_{p\in \Ob_n(\P)}
F_n(p)\rightarrow \prod_{p\in K_n} F_n(p).
$$
We also have maps
$$
\prod_{p\in \Ob_{n-1}(\P)\setminus K_{n-1}} F_{n-1}(p)
\hookrightarrow\prod_{p\in \Ob_{n-1}(\P)}
F_{n-1}(p)\stackrel{d_{n-1}}\rightarrow \prod_{p\in \Ob_n(\P)}
F_n(p)\rightarrow \prod_{p\in \Ob_n(\P)\setminus K_n} F_n(p)
$$
and
$$
\prod_{p\in K_{n-1}} F_{n-1}(p) \hookrightarrow\prod_{p\in
\Ob_{n-1}(\P)} F_{n-1}(p)\stackrel{d_{n-1}}\rightarrow \prod_{p\in
\Ob_n(\P)} F_n(p)\rightarrow \prod_{p\in K_n} F_n(p),
$$
obtained by pre and post composing the differential $d_{n-1}$ with
appropiate inclusions and projections. We denote all of them by
$d_{n-1}$. Also we fix the following
\begin{Not}
We will write elements $x\in \prod_{p\in \Ob_n(\P)} F_n(p)$ as
$x=y\oplus z$, where $y\in \prod_{p\in \Ob_n(\P)\setminus K_n}
F_n(p)$ and $z\in \prod_{p\in K_n} F_n(p)$. Also, $d_n(y\oplus
z)=(y_1\oplus z_1)\oplus(y_2+z_2)$ where $d_n(y\oplus 0)=y_1\oplus
y_2$ and $d_n(0\oplus z)=z_1\oplus z_2$.
\end{Not}
Next we introduce the Betti numbers associated to a global covering
family
\begin{Defi}\label{defi_betti_numbers_K}
Let $\P$ be a graded poset for which there exists an adequate local
covering family, and let $\K$ be a global covering family for $\P$.
Then we define, for $n\geq 0$, the $n$-th Betti number of the family
$\K$ as
$$
b^\K_0=\sum_{p\in K_0} R^p_0
$$
and
$$
b^\K_n=\sum_{p\in K_{n-1}} R^p_{n-1}-\sum_{p\in \Ob_{n-1}(\P)}
R^p_{n-1}+\sum_{p\in K_n} R^p_n
$$
for $n\geq 1$.
\end{Defi}

\begin{Rmk}
It turns out that a local covering family gives, for each
subcategory $(p_0\downarrow \P)$, a global covering family $\K$ for
$(p_0\downarrow \P)$ with $b^\K_n=0$ for $n\geq 1$ and $b^\K_0=1$.
More precisely, let $\P$ be a graded poset and let $\J$ be an
adequate local covering family for $\P$. Fix the graded poset
$\P'=(p_0\downarrow \P)$ for some object $p_0$ of $\P$. Define
$K_n=J^{p_0}_n \subseteq \Ob_n(\P')$ for $0\leq n\leq deg(p_0)$ and
empty otherwise. Then, using that $p_0$ is an initial object in
$\P'$ and Remark \ref{rmk_section}, it is straightforward that
$\K=\{K_\cdot\}_{n\geq 0}$ is a global covering family for $\P'$. In
fact, all the the maps in Definition \ref{defi_global covering
family} become isomorphisms. The integral equalities in Definition
\ref{defi_adequate_covering_family} (adequateness) correspond
exactly to the statements $b^\K_0=1$ and $b^K_n=0$ for $n\geq 1$.
\end{Rmk}
Recall once again the exact sequence
(\ref{eqn_long_exact_sequence_H^n}). From there we have that $b_n$,
the $n$-th Betti number of $\P$,  equals
$$
b_n=\dim \liminv F_{n-1}-\sum_{p\in \Ob_{n-1}(\P)} R^p_{n-1}+\dim
\liminv F_n.
$$
Because the map $\liminv F_n\rightarrow \prod_{p\in K_n} F_n(p)$ is
injective then we have
$$
\dim \liminv F_n\leq \dim \prod_{p\in K_n}
F_n(p)=\sum_{p\in K_n} R^p_n
$$
Thus, by Definition \ref{defi_betti_numbers_K},
$$
b_n\leq b^\K_n
$$
for any $n\geq 0$. An easy induction argument proves

\begin{Prop}
Let $\P$ be a poset for which exists an adequate local covering
family and a global covering family $\K$. Then
$$
b_n\leq b^\K_n
$$
and
$$
b_n-b_{n-1}+\ldots+(-1)^n b_0\leq b^\K_n-b^\K_{n-1}+\ldots+(-1)^n
b^\K_0
$$
for any $n\geq 0$.
\end{Prop}

Now we focus on the main theorem of this section, which states that
we can compute the cohomology of $\P$ using a cochain complex which
in degree $n$ has a free abelian group of rank $b^\K_n$:

\begin{Thm}\label{proposition_cohomology_cochain_global_covering}
Let $\P$ be a graded poset for which exists an adequate local
covering family and a global covering family $\K$. Then there exists
a cochain complex
$$
\xymatrix{0\ar[r]& B^\K_0 \ar[r]^{\Omega_0} & B^\K_1
\ar[r]^{\Omega_1}& B^\K_2 \ar[r]^{\Omega_2}& \ldots }
$$
with $B^\K_n\cong \Z^{b^{\K}_n}$ for $n\geq 0$, of which cohomology
is $H^*(\P;\Z)$.
\end{Thm}

Before proving the propositon we prove the following
definition-lemma

\begin{Lem}\label{lemma_betti_ses}
Let $\P$ be a graded poset for which there exists an adequate local
covering family and a global covering family $\K$. Fix $n\geq 1$,
then there is a split short exact sequence
$$
0\rightarrow\prod_{p\in \Ob_{n-1}(\P)\setminus K_{n-1}}
F_{n-1}(p)\rightarrow \prod_{p\in K_n} F_n(p)\rightarrow
B^\K_n\rightarrow 0,
$$
where $B^\K_n\cong \Z^{b^\K_n}$.
\end{Lem}
\begin{proof}
The map $\lambda:\prod_{p\in \Ob_{n-1}(\P)\setminus K_{n-1}}
F_{n-1}(p)\rightarrow \prod_{p\in K_n} F_n(p)$ is pure by definition
of global covering family, and thus its cokernel $B^\K_n$ is free
and we have a section. So, if we show that this map is a
monomorphism we obtain that the sequence is exact and the
appropriate rank of the cokernel by the definition of the number
$b^\K_n$. Take $x\in \prod_{p\in \Ob_{n-1}(\P)\setminus K_{n-1}}
F_{n-1}(p)$ with $\lambda(x)=0$. Recall the sequence
$$
\xymatrix{{\liminv F_{n-1}}\ar[r]^(.3){\iota_{n-1}}&\prod_{p\in
\Ob_{n-1}(\P)} F_{n-1}(p)\ar[r]^(.7){\omega_n}&{\liminv
F_n}\ar[r]^(.3){\iota_n}&\prod_{p\in \Ob_n(\P)} F_n(p)},
$$
which is exact at $\prod_{p\in \Ob_{n-1}(\P)} F_{n-1}(p)$. Then
$\iota_n\circ \omega_n(x\oplus 0)_p=0$ for each $p\in K_n$ and, by
Definition \ref{defi_global covering family}
(\ref{def_globalcovfam_mono}), we obtain that $\omega_n(x\oplus
0)=0$. By exactness there exist $y\in \liminv F_{n-1}$ with
$\iota_{n-1}(y)=x\oplus0$. But, by definition, $(x\oplus0)_p=0$ for
each $p\in K_{n-1}$. Then, by Definition \ref{defi_global covering
family} (\ref{def_globalcovfam_mono}) again, we obtain that $y=0$
and $x=0$.
\end{proof}

The leftmost square in the following diagam commutes, and thus, we
can find an arrow which closes the rightmost square. We denote this
arrow by $\Omega_n$, and it is the differential on $B_\cdot$ induced
by the differential $d_n$:
$$
\xymatrix{ 0\ar[r]&\prod_{p\in \Ob_{n-1}(\P)\setminus K_{n-1}}
F_{n-1}(p)\ar[r]^{d_{n-1}}\ar[d]^{-d_{n-1}}& \prod_{p\in K_n}
F_n(p)\ar[r]\ar[d]^{d_n}&
B^\K_n\ar[r]\ar@{-->}[d]^{\Omega_n}&0\\
0\ar[r]&\prod_{p\in \Ob_n(\P)\setminus K_n} F_n(p)\ar[r]^{d_n}&
\prod_{p\in K_{n+1}} F_{n+1}(p)\ar[r]& B^\K_{n+1}\ar[r]&0.}
$$
The maps $\Omega_\cdot$ verify $\Omega_{n+1}\circ \Omega_n=0$ for
each $n\geq 0$, as can be seen by using the preceding diagram. We
have then a chain complex
$$
\xymatrix{0\ar[r]& B^\K_0 \ar[r]^{\Omega_0} & B^\K_1
\ar[r]^{\Omega_1}& B^\K_2 \ar[r]^{\Omega_2}& \ldots },
$$
and we claim that the homology of this chain complex is the same as
that of the chain complex in Theorem
\ref{proposition_cohomology_cochain_covering}.

Thus to prove the proposition we shall build for each $n\geq 1$ an
epimorphism
$$
\psi:\ker\Omega_n\rightarrow \ker d_n/\im d_{n-1}
$$
of which kernel is $\im \Omega_{n-1}$. We denote the class of $x\in
\prod_{p\in K_n} F_n(p)$ in $B^\K_n$ as $\overline{x}$, and the
class of $x\in \ker d_n$ in $\ker d_n/\im d_{n-1}$ by
$\overline{x}$.

Thus, take $\overline{x}\in B^\K_n$ such that
$\Omega_n(\overline{x})=0$. Then there is $y\in \prod_{p\in
\Ob_n(\P)\setminus K_n} F_n(p)$ such that $d_n(x)=d_n(y)$. If
$d_n(y\oplus 0)=y_1\oplus y_2$ and $d_n(0\oplus x)=x_1\oplus x_2$
then $y_2=x_2$. Define
$$
\psi(\overline{x})=\overline{-y\oplus x}.
$$
Then $d_n(-y\oplus x)=-(y_1\oplus y_2)+(x_1\oplus x_2)=x_1-y_1\oplus
0$. Therefore $\iota_{n+1}(\omega_{n+1}(-y\oplus x))=x_1-y_1\oplus
0$ and, as the map
$$
\liminv F_{n+1}\stackrel{\iota_{n+1}}\rightarrow \prod_{p\in
\Ob_{n+1}(\P)} F_{n+1}(p)\rightarrow \prod_{p\in \K_{n+1}(\P)}
F_{n+1}(p)
$$
is injective by definition of global covering family, we obtain that
$x_1=y_1$ and $d_n(-y\oplus x)=0$, i.e., $-y\oplus x\in \ker d_n$.
Notice that the element $y$ chosen above is unique: if $y'\in
\prod_{p\in \Ob_n(\P)\setminus K_n} F_n(p)$ and $d_n(x)=d_n(y)$ then
$d_n(y-y')=0$ and by Lemma \ref{lemma_betti_ses} $y=y'$.

Now we prove that $\psi(\overline{x})$ does not depend on the chosen
representative $x$. Thus take $x'=x+z_2$, where $z\in \prod_{p\in
\Ob_{n-1}(\P)\setminus K_{n-1}} F_{n-1}(p)$ and $d_{n-1}(0\oplus
z)=z_1\oplus z_2$. Then $0=d_n(d_{n-1}(0\oplus
z))=z_{1,1}+z_{2,1}\oplus z_{1,2}+z_{2,2}$. Moreover,
$d_n(y-z_1\oplus 0)=y_1-z_{1,1}\oplus y_2-z_{1,2}$ and thus
$d_n(y-z_1)=y_2-z_{1,2}=y_2+z_{2,2}=x_2+z_{2,2}=d_n(x')$. Then,
$$
\psi(\overline{x'})=\overline{-y+z_1\oplus x'}.
$$
But notice that $d_{n-1}(z\oplus 0)=z_1\oplus z_2=(-y+z_1\oplus x')
- (-y\oplus x)$ and thus $\overline{-y+z_1\oplus
x'}=\overline{-y\oplus x}$ in $\ker d_n/ \im d_{n-1}$.

It is clear that $\psi$ is a homomorphism. Next we prove that it is
an epimorphism. Take $\overline{a\oplus b}\in \ker d_n/ \im
d_{n-1}$. Then $0=d_n(a\oplus b)=a_1+b_1\oplus a_2+b_2$. Consider
$b\in \prod_{p\in K_n} F_n(p)$. Notice that $d_n(a\oplus
0)=a_1\oplus a_2$ and $d_n(0\oplus b)=b_1\oplus b_2$. Hence
$d_n(b)=b_2=-a_2=d_n(-a)$, and thus $\Omega_n(\overline{b})=0$ and
$\psi(\overline{b})=\overline{a\oplus b}$.

To finish the proof of the proposition we show that $\ker \psi=\im
\Omega_{n-1}$. First take $\overline{x}\in B^\K_n$ with
$\overline{x}=\Omega_{n-1}(\overline{y})$ and $\overline{y}\in
B^\K_{n-1}$. There is $z\in \prod_{p\in \Ob_{n-1}(\P)\setminus
K_{n-1}} F_{n-1}(p)$ with $d_{n-1}(z)=x-d_{n-1}(y)$. If
$d_{n-1}(z\oplus 0)=z_1\oplus z_2$ and $d_{n-1}(0\oplus y)=y_1\oplus
y_2$ then $z_2=x-y_2$ and $d_n(x)=x_2=z_{2,2}+y_{2,2}$. Take
$t=-z_1-y_1\in \prod_{p\in \Ob_{n-1}(\P)\setminus K_{n-1}}
F_{n-1}(p)$. Therefore $d_n(t,0)=-z_{1,1}-y_{1,1}\oplus
-z_{1,2}-y_{1,2}$, $d_n(t)=-z_{1,2}-y_{1,2}=z_{2,2}+y_{2,2}=x_2$ and
$\psi(\overline{x})=\overline{-t\oplus x}$ with $d_n(z\oplus
y)=z_1+y_1\oplus z_2+y_2=-t\oplus x$. This shows that $\im
\Omega_{n-1}\subseteq \ker \psi$.

Finally we prove that $\ker \psi\subseteq \im \Omega_{n-1}$. Take
$\overline{x}\in B^\K_n$ such that $\psi(\overline{x})=0$ in $\ker
d_n/ \im d_{n-1}$. This means that if $y\in \prod_{p\in
\Ob_n(\P)\setminus K_n} F_n(p)$ is such that $d_n(x)=d_n(y)$ then
there is $a\oplus b\in \prod_{p\in \Ob_{n-1}(\P)} F_{n-1}(p)$ with
$d_{n-1}(a\oplus b)=-y\oplus x$. Hence $a\in \prod_{p\in
\Ob_{n-1}(\P)\setminus K_{n-1}} F_{n-1}(p)$ and $b\in \prod_{p\in
\K_{n-1}(\P)} F_{n-1}(p)$ are such that $d_{n-1}(b)=x-d_{n-1}(a)$
and thus $\Omega_{n-1}(\overline{b})=\overline{x}$.

%% file: thing-like.tex
\section{Simplex-like posets}\label{section_locally delta}
In this section we will show that for a big family of posets, which
includes simplicial complexes, there exists local covering families.
We start defining posets which locally are like simplicial
complexes:
\begin{Defi}\label{defi_simplex n}
The category $\triangle_n$ has as objects the faces of the standard
$n$-dimensional simplex and arrows the inclusions among faces. For
example $\triangle_2$ is the category:
$$
\xymatrix{ &-\ar@{<-}[r]\ar@{<-}[rd]&\cdot\\
\triangle\ar@{<-}[r]\ar@{<-}[ru]\ar@{<-}[rd]&-\ar@{<-}[ru]\ar@{<-}[rd]&\cdot\\
&-\ar@{<-}[r]\ar@{<-}[ru]&\cdot.}
$$
\end{Defi}

%\begin{Defi}\label{defi_cube n}
%The category $\Box_n$ has as objects the faces of the standard
%$n$-dimensional cube $[0,1]^n$ and arrows the inclusions among faces. For
%example $\Box_2$ is the category:
%$$
%\xymatrix{ &-\ar@{<-}[r]\ar@{<-}[rd]&\cdot\\
%\Box\ar@{<-}[r]\ar@{<-}[ru]\ar@{<-}[rd]\ar@{<-}[rdd]&-\ar@{<-}[r]\ar@{<-}[rd]&\cdot\\
%&-\ar@{<-}[r]\ar@{<-}[rd]&\cdot\\
%&-\ar@{<-}[r]\ar@{<-}[ruuu]&\cdot.}
%$$
%\end{Defi}

\begin{Defi}\label{simplex-like poset}
Let $\P$ be a poset. Then $\P$ is \emph{simplex-like} if for all
$p\in \Ob(\P)$ the category $(\P\downarrow p)$ is isomorphic to
$\triangle_n$ for some $n$.
\end{Defi}

%\begin{Defi}\label{cubical-like poset}
%Let $\P$ be a poset. Then $\P$ is \emph{cubical-like} if for all
%$p\in \Ob(\P)$ the category $(\P\downarrow p)$ is isomorphic to
%$\Box_n$ for some $n$.
%\end{Defi}

Of course simplicial complexes are simplex-like posets. In fact
(\cite[3.1]{garrett}), simplicial complexes are exactly the
simplex-like posets such that any two elements which have a lower
bound have an infimum, i.e., a greatest lower bound. Another
examples of simplex-like posets are barycentric subdivisions and, in
general, subdivision categories in the sense of \cite{markus}. We
have:

\begin{Lem}\label{thinglike_is_graded}
Let $\P$ be a simplex-like poset. Then $\P^{op}$ is a graded poset.
\end{Lem}
\begin{proof}
To see that $\P$ is graded recall that for any $p\in \Ob(\P)$ the
subcategory $(\P\downarrow p)$ is isomorphic to $\triangle_n$ for
some $n\geq 0$. Define $deg(p)=n$. Then $deg:\Ob(\P^{op})\rightarrow
\Z$ is a degree function which assigns $0$ to maximal elements and
$\P$ is graded.
\end{proof}

Next, we will define local covering families for $\triangle_n^{op}$.
Due to the local shape of simplex-like posets this local covering
family for $\triangle_n^{op}$ will allow us to define a local
covering family for the opposite category of a simplex-like poset.

\begin{Lem}\label{lem_coveringfamily simplex}
There exists an adequate local covering family for
$\triangle_n^{op}$ for each $n\geq 0$.
\end{Lem}
\begin{proof}
Fix a total order $v_0<v_1<\ldots<v_n$ for the vertices of
$\triangle_n$. For each face $\sigma$ of $\triangle_n$ and each
$0\leq m\leq deg(\sigma)$ we must define a subset $\J^\sigma_m$ of
the $m$-dimensional faces of $\sigma$. Take the greatest vertex
contained in $\sigma$ and then define $\J^\sigma_m$ as those
$m$-dimensional faces of $\sigma$ which contain this vertex. Then it
is straightforward that $\J$ is an adequate local covering family of
$\triangle_n^{op}$, and that
$$
R^{\sigma}_m=\sum_{l=0}^m (-1)^{(m-l)}\big(\!\begin{smallmatrix}
deg(\sigma)+1
\\ l\end{smallmatrix}\!\big)
$$
for $\sigma\in \triangle_n$ and $0\leq m\leq deg(\sigma)$.
\end{proof}

Now we reach the main result of this section:

\begin{Lem}\label{lem_coveringfamily_for_locallyDelta}
If $\P$ is a simplex-like poset then there is an adequate local
covering family for the graded poset $\P^{op}$.
\end{Lem}
\begin{proof}
By definition there are isomorphisms of categories $(P\downarrow
p)^{op}\simeq (p\downarrow \P^{op})\simeq \triangle_n^{op}$ for each
$p\in \Ob(\P)$, and we know that $\triangle_n^{op}$ can be equipped
with an adequate local covering family. To build an adequate local
covering family $$\K=\{K^{i_0}_p\}_{i_0\in \Ob(\P^{op})\text{,
$0\leq p\leq deg(i_0)$}}$$we just have to choose appropriately the
isomorphisms $(P\downarrow p)\simeq \triangle_n$.

Consider the degree function $deg:\Ob(\P^{op})\rightarrow \Z$
defined in Lemma \ref{thinglike_is_graded} and the set $T$ of the
maximal elements of $\P^{op}$, i.e., $T=\{p\in
\Ob(\P^{op})|deg(p)=0\}$. Choose a total order $<$ for $T$ (suppose
$T$ is finite or use the Axiom of Choice \cite{jech}). Then, given
$p\in \Ob(\P^{op})$, consider the subset $(p\downarrow
\P^{op})_0\subseteq T$ and the restriction $((p\downarrow
\P^{op})_0,<)$ of the total order from $T$. There is a unique
isomorphism
$$
\varphi_p:(p\downarrow \P^{op})\simeq \triangle_{deg(p)}^{op}
$$
which induces an order preserving map
$$
((p\downarrow \P^{op})_0,<)\simeq
(\triangle_{deg(p)})_0=\{v_0,v_1,v_2,\ldots,v_{deg(p)}\}.
$$

Denote by $\J$ the local covering family for
$\triangle_{deg(p)}^{op}$ of Lemma \ref{lem_coveringfamily simplex}
and define, for $0\leq m\leq deg(p)$,
$$
K^{p}_m=\varphi_{p}^{-1}(J^{\varphi(p)}_m).
$$

Then $\K$ fulfills condition \ref{defi_covering_a}) of Definition
\ref{defi_covering_family} because for $0\leq m<deg(p)$
$$
(p\downarrow
\P^{op})_m=\varphi_{p}^{-1}((\triangle_{deg(p)}^{op})_m)=\varphi_{p}^{-1}(\bigcup_{i\in
J^{\varphi_{p}(p)}_{m+1}} (i\downarrow \triangle_{deg(p)}^{op})_m)=
$$
$$=\bigcup_{\varphi_{p}^{-1}(i)\in \varphi_{p}^{-1}(K^{p}_{m+1})}
\varphi_{p}^{-1}((i\downarrow
\triangle_{deg(p)}^{op})_m)=\bigcup_{i\in K^{p}_{m+1}} (i\downarrow
\P^{op})_m.
$$

To check condition \ref{defi_covering_b}) of Definition
\ref{defi_covering_family} take $q\in K^{p}_{m+1}$ for some $0\leq
m<deg(p)$ and call $\J'$ to the local covering family for
$\triangle_{m+1}^{op}$ of Lemma \ref{lem_coveringfamily simplex}. We
want to see that $K^q_m\subseteq K^{p}_m$. Recalling the natural
inclusion $(q\downarrow \P^{op})\subseteq (p\downarrow \P^{op})$
this is equivalent to
$$
\varphi_q^{-1}(J'^{\varphi_q(q)}_m)\subseteq
\varphi_{p}^{-1}(J^{\varphi_{p}(p)}_m)
$$
and to
$$
\psi(J'^{\varphi_q(q)}_m)\subseteq J^{\varphi_{p}(p)}_m
$$
where $\psi=\varphi_{p}\circ \varphi_q^{-1}$. By construction $\psi$
is order preserving and thus this inclusion holds.
\end{proof}

Notice that comparing the general expression for $R^p_n$ in the
proof of Lemma \ref{lem_coveringfamily simplex} with the binomial
expansion of $(1-1)^{deg(p)+1}=0$ we obtain that $R^p_{deg(p)}=1$
for any object $p$ in a simplex-like poset. The following are
re-statements of results about local covering families applied to
simplex-like posets:

\begin{Lem}\label{lema_adequate_covering_family_simplex-like posets}
Let $\P$ be a simplex-like poset and consider the graded poset
$\P^{op}$ (for which exists an adequate local covering family by
Lemma {\rm\ref{lem_coveringfamily_for_locallyDelta}}). Let $\K$ be a
global covering family for $\P^{op}$. Then the Betti numbers of the
family $\K$ are given by
$$
b^\K_0=|K_0|,
$$
and for $n\geq 1$, by
$$
b^\K_n=|K_{n-1}|-|\Ob_{n-1}(\P)|+|K_n|.
$$
\end{Lem}

\begin{Lem}\label{lemma_analogue eulerchar_simplex-like posets}
Let $\P$ a simplex-like poset. Then
$$
\sum_i (-1)^i\dim H^i(\P;\Z)=\sum_i (-1)^i |\Ob_i(\P)|.
$$
\end{Lem}

Again because $R^p_{deg(p)}=1$, Proposition \ref{cohomology with
adequate local covering family} gives $\dim F_{deg(p)}(p)=1$ and
threfore $F_{deg(p)}(p)=\Z$. We may identify the object $p\in \P_n$
with its image $v_0v_1\ldots v_n$ under the isomorphism
$(\P\downarrow p)\cong\triangle_n$ built in Lemma
\ref{lema_adequate_covering_family_simplex-like posets}. Then the
expression for the differential in Theorem
\ref{proposition_cohomology_cochain_covering} becomes the familiar
expression for  simplicial complexes:

\begin{Prop}\label{proposition_cohomology_cochain_covering_simplex-like}
Let $\P$ be a simplex-like poset. Then there exists a co-chain
complex
$$
0\rightarrow \Z^{\Ob_0(\P)}\stackrel{d_0}\rightarrow \Z^{\Ob_1(\P)}
\stackrel{d_1}\rightarrow \Z^{\Ob_2(\P)} \stackrel{d_2}\rightarrow
\ldots
$$
of which cohomology is $H^*(\P;\Z)$. Under the identifications
described above we have for $p=v_0v_1\ldots v_n\in \Ob_n(\P)$,
$n\geq 1$ and $x\in \Z^{\Ob_{n-1}(\P)}$
$$
d_{n-1}(x)_p=\sum_{j=0}^n (-1)^{n-j} x_{v_0v_1\ldots\hat{v_j}\ldots
v_n}
$$
\end{Prop}
\begin{proof}
Take $p$ and $x$ as in the statement. The short exact sequence
considered in Remark \ref{rmk_dimensions}
$$
0\rightarrow F_{n-1}(p)\stackrel{\lambda_{p}}\rightarrow
\ker_{F_{n-1}}'(p)\stackrel{\pi_{p}}\rightarrow F_n(p)\rightarrow 0
$$
takes the form
$$
0\rightarrow F_{n-1}(v_0v_1\ldots
v_n)\stackrel{\lambda_{p}}\rightarrow \prod_{j=0}^n
F_{n-1}(v_0\ldots \hat{v_j}\ldots v_n)\stackrel{\pi_{p}}\rightarrow
F_n(v_0v_1\ldots v_n)\rightarrow 0.
$$
By the definition of the differential $d_{n-1}$ we have that
$d_{n-1}(x)_p=\pi_p(y)$, where $y=x|_{\{v_0\ldots \hat{v_j}\ldots
v_n\text{, $j=0,\ldots,n$}\}}\in \prod_{j=0}^n F_{n-1}(v_0\ldots
\hat{v_j}\ldots v_n)$. As in Remark \ref{rmk_section} we identify
$\pi_p(y)$ with its unique pre-image $y'\in \prod_{j=0}^n
F_{n-1}(v_0\ldots \hat{v_j}\ldots v_n)$ such that $y'_{v_0\ldots
\hat{v_j}\ldots v_n}=0$ for $j=0,\ldots,n-1$. Thus, if we obtain
$z\in F_{n-1}(v_0v_1\ldots v_n)$ such that $\lambda_p(z)_{v_0\ldots
\hat{v_j}\ldots v_n}=y_{v_0\ldots \hat{v_j}\ldots v_n}$ for
$j=0,\ldots,n-1$, we will have $\pi_p(y)\equiv y'=y-\lambda_p(z)$.

The rest of the proof is by induction on $n$. For $n=1$ the short
exact sequence is
$$
0\rightarrow \Z_{v_0v_1}\stackrel{\lambda_{p}}\rightarrow
\Z_{v_0}\times \Z_{v_1}\stackrel{\pi_{p}}\rightarrow
\Z_{v_0v_1}\rightarrow 0,
$$
where $\lambda_p(n)=(n,n)$. Then, if $y=(y_{v_0},y_{v_1})$, we have
$z=y_{v_1}$ and
$$y'=y-\lambda_p(z)=(y_{v_0},y_{v_1})-(y_{v_1},y_{v_1})=(y_{v_0}-y_{v_1},0).$$
Now we make the inductive step for $n\geq 2$. Recall that we want
$z$ such that $\lambda_p(z)_{v_0\ldots \hat{v_j}\ldots
v_n}=y_{v_0\ldots \hat{v_j}\ldots v_n}$ for $j=0,\ldots,n-1$. Again
by Remark \ref{rmk_section} we identify
\begin{equation}\label{equ_identificacionA}
F_{n-1}(v_0v_1\ldots v_n)\cong \prod_{j=0}^{n-1}
F_{n-2}(v_0v_1\ldots \hat{v_j}\ldots v_{n-1}),
\end{equation}
and, for $j=0,\ldots,n-1$,
\begin{equation}\label{equ_identificacionB}
F_{n-1}(v_0v_1\ldots\hat{v_j}\ldots v_n)\cong F_{n-2}(v_0v_1\ldots
\hat{v_j}\ldots v_{n-1}),
\end{equation}
and
\begin{equation}\label{equ_identificacionC}
F_{n-1}(v_0v_1\ldots v_{n-1})\cong F_{n-2}(v_0v_1\ldots v_{n-2}).
\end{equation}

Recall that the ``diagonal'' $\lambda_p$ is given by $\prod_{j=0}^n
F_{n-1}(v_0\ldots v_n\rightarrow v_0\ldots \hat{v_j}\ldots v_n)$.
For $j_0\in\{0,\ldots,n-1\}$, and with the identifications
(\ref{equ_identificacionA}) and (\ref{equ_identificacionB}), the map
$F_{n-1}(v_0\ldots v_n\rightarrow v_0\ldots \hat{v_{j_0}}\ldots
v_n)$ becomes exactly the projection
$$
\prod_{j=0}^{n-1} F_{n-2}(v_0v_1\ldots \hat{v_j}\ldots
v_{n-1})\rightarrow F_{n-2}(v_0v_1\ldots \hat{v_{j_0}}\ldots
v_{n-1}).
$$
Thus, $z\in F_{n-1}(v_0v_1\ldots v_n)$ is such that $z_{v_0v_1\ldots
\hat{v_j}\ldots v_{n-1}}=y_{v_0v_1\ldots \hat{v_j}\ldots v_n}$
(again under the identifications (\ref{equ_identificacionA}) and
(\ref{equ_identificacionB})). The only thing left to compute is
$F_{n-1}(v_0\ldots v_n\rightarrow v_0v_1\ldots v_{n-1})(z)$ under
the identifications (\ref{equ_identificacionA}) and
(\ref{equ_identificacionC}). Define $q=v_0v_1\ldots v_{n-1}$ and
consider the short exact sequence given by Remark
\ref{rmk_dimensions}
$$
0\rightarrow F_{n-2}(v_0v_1\ldots
v_{n-1})\stackrel{\lambda_q}\rightarrow \prod_{j=0}^{n-1}
F_{n-2}(v_0v_1\ldots \hat{v_j}\ldots
v_{n-1})\stackrel{\pi_q}\rightarrow F_{n-1}(v_0v_1\ldots
v_{n-1})\rightarrow 0.
$$
Now, by the induction hypothesis, $\pi_q(z)=\sum_{j=0}^{n-1}
(-1)^{n-1-j} z_{v_0v_1\ldots\hat{v_j}\ldots v_{n-1}}$. Thus,
finally,
$$
y'=y-z=(y_{v_0\ldots v_{n-1}}-\sum_{j=0}^{n-1} (-1)^{n-1-j}
z_{v_0v_1\ldots\hat{v_j}\ldots v_{n-1}},0,\ldots,0)
$$
which equals
$$
y'=(\sum_{j=0}^n (-1)^{n-j} y_{v_0v_1\ldots\hat{v_j}\ldots
v_n},0,\ldots,0).
$$
\end{proof}
Before finishing this section we show how the symmetry of the
category $\triangle_n$ is useful when checking if a given family of
subsets of objects $\K=\{K_n\}_{n\geq 0}$ fulfills the properties of
a global covering family. Consider a simplex like poset $\P$ and
call $\J$ to the local covering family for $(p\downarrow
\P^{op})\simeq \triangle_n^{op}$ built in Lemma
\ref{lem_coveringfamily_for_locallyDelta}. Denote $p=v_0v_1\ldots
v_n$ and assume $deg(p)=n\geq 1$, then there are isomorphisms
(Remark \ref{rmk_section})
$$
F_{n-1}(v_0v_1\ldots v_n)\stackrel{\cong}\rightarrow
\prod_{j=0}^{n-1} F_{n-1}(v_0\ldots \hat{v_j}\ldots v_n)
$$
and
$$
F_n(v_0v_1\ldots v_n)\stackrel{\cong}\rightarrow
F_{n-1}(v_0v_1\ldots v_{n-1}).
$$
By letting the symmetric group $\Sigma_{n+1}$ act on
$\triangle_n^{op}$ we map $\J$ to other local covering families for
$(p\downarrow \P^{op})$. In particular we can permute $v_n$ with any
of the objects $\{v_0,\ldots ,v_{n-1}\}$. Recall that the functors
$F_\cdot$ do not depend on the local covering family chosen. We may
summarize this as

\begin{Rmk}\label{Remark_symmetry_triangle_n}
Let $\P$ be a simplex-like poset and let $p=v_0v_1\ldots v_n$ be an
object with $deg(p)=n\geq 1$. The maps
$$
F_{n-1}(v_0v_1\ldots v_n)\rightarrow \prod_{j=0,j\neq j_0}^{n}
F_{n-1}(v_0\ldots \hat{v_j}\ldots v_n)
$$
and
$$
F_n(v_0v_1\ldots v_n)\rightarrow
F_{n-1}(v_0\ldots\hat{v_{j_0}}\ldots v_{n-1}).
$$
are both isomorphisms for $j_0=0,\ldots,n$ .
\end{Rmk}

\begin{Ex}\label{ex_projective_plane}
In this example we consider the real projective plane $\R P^2$ and a
poset model $\P$ of it. It has four $2$-cells, six edges and three
vertices:

\begin{tabular}{ccccccccccc}
{\entrymodifiers={}
$$
\xymatrix@=50pt{\ar@{}[rd]^<w\ar@{-}[r]\ar@{-}[dr]^(.4)d&\ar@{}[d]^B\ar@{<-}[r]^b&\ar[d]\ar@{-}[ld]^(.4)e\ar@{}[ld]^<v\\
\ar@{}[r]^A\ar@{-}[u]^a&\ar@{}[r]^<x&\ar@{-}[d]^a\ar@{}[l]^C\\
\ar@{}[ru]^<v\ar@{-}[ru]^(.4)c\ar[u]\ar[r]_b&\ar@{}^D[u]\ar@{-}[r]&\ar@{-}[lu]^(.4)f\ar@{}[lu]^<w}
$$}
&&&&&&&&&&
{\tiny{
$$
\xymatrix@R=10pt@C=50pt{A \ar@{<-}[r]\ar@{<-}[rdd]\ar@{<-}[rddd]& a\ar@{<-}[r]\ar@{<-}[rd]&v\\
B\ar@{<-}[r]\ar@{<-}[rdd]\ar@{<-}[rddd]&b\ar@{<-}[ru]\ar@{<-}[r]&w\\
&c\ar@{<-}[ruu]\ar@{<-}[rddd]&\\
&d\ar@{<-}[ruu]\ar@{<-}[rdd]&\\
C\ar@{<-}[ruuuu]\ar@{<-}[r]\ar@{<-}[rd]&e\ar@{<-}[ruuuu]\ar@{<-}[rd]&\\
D\ar@{<-}[ruuuu]\ar@{<-}[ruuu]\ar@{<-}[r]&f\ar@{<-}[ruuuu]\ar@{<-}[r]&x}
$$}}
\end{tabular}

This poset is a simplex-like poset and thus there is a local
covering family for its opposite category $\P^{op}$. Notice that it
is not a simplicial complex as, for example, there are two edges
($c$ and $e$) with common vertices $v$ and $x$. Applying Proposition
\ref{proposition_cohomology_cochain_covering_simplex-like} we obtain
that its cohomology is computed by a co-chain complex:

$$
0\rightarrow
\Z^3\stackrel{d_0}\rightarrow\Z^6\stackrel{d_1}\rightarrow
\Z^4\rightarrow 0.
$$

Proposition
\ref{proposition_cohomology_cochain_covering_simplex-like} gives a
description of the differentials $d_0$ and $d_1$. Now consider the
following family $\K=\{K_0,K_1,K_2\}$ of subsets of objects of $\P$:
\begin{eqnarray}
&K_0=\{x\}, \nonumber\\
&K_1=\{d,e,f\},\nonumber\\
&K_2=\{A,B,C,D\}.\nonumber
\end{eqnarray}
Using Remark \ref{Remark_symmetry_triangle_n} it is straightforward
that $\K$ is a global covering family. Its Betti numbers are
$b^\K_0=1$, $b^\K_1=1$ and $b^\K_2=1$ and thus, by Theorem
\ref{proposition_cohomology_cochain_global_covering}, the cohomology
of $\R P^2$ is that of a cochain complex
$$
0\rightarrow
\Z\stackrel{\Omega_0}\rightarrow\Z\stackrel{\Omega_1}\rightarrow
\Z\rightarrow 0.
$$
The induced differentials are easily computed to be $\Omega_0\equiv
0$ and $\Omega_1\equiv\times 2$.
\end{Ex}

%% file: morse.tex
\section{Morse theory}\label{section_morse}
In this section we show how any Morse function on a simplicial
complex gives rise to a global covering family. The setup for this
section is the discrete Morse theory for $CW$-complexes that Forman
introduces in \cite{forman}. We will restrict ourselves to
simplicial complexes, as the same author does in the user's guide
\cite{forman_survey}.

We start introducing the concept of Morse function in this context.
Suppose $\P$ is a given simplicial complex (in this section all
simplicial complexes are assumed to have a finite number of
vertices). Also, for a simplex $p$ we write $p^n$ to denote that $n$
is the dimension of $p$. A function $f:\P\rightarrow \R$ is called a
\emph{Morse} function if for each $p^n\in \P$
$$
|\{q^{n+1}>p|f(q)\leq f(p) \}|\leq 1
$$
and
$$
|\{q^{n-1}<p|f(p)\leq f(q) \}|\leq 1.
$$
This means roughly that $f$ increases as the dimension of the
simplices increase, with at most one exception, locally, at each
simplex. We reproduce here the basic result \cite[Lemma 2.5]{forman}
\begin{Lem}\label{lemma_morse_disjuntive}
Let $\P$ be a simplicial complex and $f:\P\rightarrow \R$ a Morse
function. Then for any simplex $p^n$ either
$$
|\{q^{n+1}>p| f(q)\leq f(p) \}|=0
$$
or
$$
|\{q^{n-1}<p| f(p)\leq f(q) \}|=0.
$$
\end{Lem}

If both conditions in the lemma hold for a simplex then we call it
critical, i.e., for $p^n$ in $\P$ we say that $p$ is \emph{critical}
if
$$
|\{q^{n+1}>p|f(q)\leq f(p) \}|=0
$$
and
$$
|\{q^{n-1}<p|f(p)\leq f(q) \}|=0.
$$

Now, we come back to the setup of posets. Because $\P$ is a
simplicial complex then it is a simplex-like poset, and by Lemma
\ref{lem_coveringfamily_for_locallyDelta} there is an adequate local
covering family for $\P^{op}$. In the next proposition we see how
any Morse function $f:\P\rightarrow \R$ gives rise to a global
covering family on $\P^{op}$.

\begin{Prop}\label{morse_prop_modelo}
Let $\P$ be a simplicial complex and $f:\P\rightarrow \R$ a Morse
function. Then there is a global covering family $\K$ for $\P^{op}$.
Moreover, for each $n\geq 0$,
$$
b^\K_n=|\{\text{critical simplices of dimension $n$}\}|.
$$
\end{Prop}
\begin{proof}
According to Lemma \ref{lemma_morse_disjuntive} we can divide the
simplices of dimension $n\geq 0$, $\P_n$, in the following disjoint
sets
$$
C_n=\{p^n||\{q^{n+1}>p|f(q)\leq f(p) \}|=0\text{ and }
|\{q^{n-1}<p|f(p)\leq f(q) \}|=0\},
$$
$$
D_n=\{p^n||\{q^{n-1}<p|f(p)\leq f(q) \}|=1\},
$$
and
$$
E_n=\{p^n||\{q^{n+1}>p|f(q)\leq f(p) \}|=1\}.
$$
The set $C_n$ consists of the $n$-dimensional critical simplices.
For $n\geq 1$ there is a bijection
$$
E_{n-1}\rightarrow D_n
$$
which maps $p^{n-1}\in E_{n-1}$ to the unique $n$-simplex $q^n$ with
$f(q)\leq f(p)$. Now define $K_n=C_n\cup D_n$ (disjoint union) for
each $n\geq 0$. Notice that if $\K=\{K_n\}_{n\geq 0}$ were a global
covering family then
$$
b^\K_0=|K_0|=|C_0|+|D_0|=|C_0|
$$
as $D_0=\emptyset$ and, for $n\geq 1$,
$$
b^\K_n=|K_{n-1}|-|\P_{n-1}|+|K_n|.
$$
Because of the bijection $E_{n-1}\cong D_n$ this equals
$$
b^\K_n=|C_{n-1}|+|D_{n-1}|-(|C_{n-1}|+|D_{n-1}|+|E_{n-1}|)+|C_n|+|D_n|=|C_n|.
$$
Fix $n\geq 0$. We show that the restriction map
$$
w:\liminv F_n\rightarrow \prod_{p\in K_n} F_n(p)
$$
is a monomorphism, where $F_n:\P^{op}\rightarrow \Ab$ are the
functors obtained from Lemma \ref{cohomology with adequate local
covering family} applied to $\P^{op}$.

Take $\psi\in \liminv F_n=\hom_{\Ab^{\P}}(c_{\Z},F_n)$ which goes to
zero by the restriction map $w$. If we prove that $\psi(q)=0$ for
each simplex of dimension $n+1$ then $\psi(p)=0$ for every simplex
of dimension $n$. To see this consider any simplex $p^n$. If $p$ is
not the face of any $(n+1)$-simplex then $p\in K_n$ and $\psi(p)=0$
(as $\psi$ is in the kernel of $w$). If there exists $q^{n+1}$ with
$q>p$ then $\psi(q)=0$ by hypothesis and then
$\psi(p)=\psi(q\rightarrow p)(\psi(q))=\psi(q\rightarrow p)(0)=0$.
Finally, as $F_n$ is $n$-condensed, if $\psi(p)=0$ for each
$n$-simplex $p$ then $\psi=0$.

Now we prove that $\psi(q)=0$ for any $(n+1)$-simplex $q$. Recall
that we are assuming that the set of vertices of $\P$, and thus $\P$
itself, is finite. We consider the total ordered finite set
$f(\P_{n+1})$ and make induction on it. The base case is a
$(n+1)$-simplex $q$ such that $f(q)=min\{ f(\P_{n+1})\}$. There are
$n+2$ $n$-dimensional faces of $q$. If at least $n+1$ of these faces
are in $K_n$ then $\psi(p)=0$ for each of these faces and, by
\ref{Remark_symmetry_triangle_n}, $\psi(q)=0$.

Now suppose that there exists at least two $n$-dimensional faces of
$q$ which are not in $K_n$. By the definition of Morse function one
of the values of $f$ in these faces is strictly smaller than $f(q)$.
Call this face $p$ so $f(q)>f(p)$. As $p$ belongs to $E_n$ there
exists an $(n+1)$-dimensional simplex $q'$ such that $q'>p$ and
$f(q')\leq f(p)$. But then we obtain $f(q)>f(p)\geq f(q')$. This
contradicts the fact that $f(q)$ is minimun among $(n+1)$-simplices
and thus it has to be the case that there are at least $n+1$
$n$-dimensional faces of $q$ which are in $K_n$.

Next we do the induction step: take an $(n+1)$-dimensional simplex
$q$. By definition of Morse function at least $n+1$ of the $n+2$
$n$-dimensional faces $p$ of $q$ satisfy $f(q)>f(p)$. For each one
of these $n+1$ faces $p$ either $p\in K_n$ and $\psi(p)=0$ or $p\in
E_n$ and there exists a $(n+1)$-dimensional simplex $q'$ with
$f(q')\leq f(p)$. In this last case we obtain $f(q)>f(p)\geq f(q')$
and then by the induction hypothesis $\psi(q')=0$ and
$\psi(p)=\psi(q'\rightarrow p)(\psi(q'))=0$. Thus, $\psi$ is zero in
at least $n+1$ $n$-dimensional faces of $q$ and by
\ref{Remark_symmetry_triangle_n} $\psi(q)=0$.

Now we prove that for $n\geq 1$ the map
$$
w:\prod_{p\in \Ob_{n-1}(\P)\setminus K_{n-1}} F_{n-1}(p)\rightarrow
\prod_{p\in K_n} F_n(p)
$$
is pure. Take $y\in \prod_{p\in K_n} F_n(p)$ , $m\geq 1$ and $z\in
\prod_{p\in \Ob_{n-1}(\P)\setminus K_{n-1}} F_{n-1}(p)$ with $m\cdot
y=w(z)$. We want to prove that $m|z_p$ for each $p\in
\Ob_{n-1}(\P)\setminus K_{n-1}$. We do this by induction on
$f(\Ob_{n-1}(\P)\setminus K_{n-1})$. The base case is a
$(n-1)$-simplex $p$ not in $K_{n-1}$ such that $f(p)$ is a minimum
value. As $p$ is in $E_{n-1}$ there is $q^n\in D_n\subset K_n$ with
$f(q)\leq f(p)$. For any of the $n$ faces $p'$ of $q$ different from
$p$ we have $f(q)>f(p')$ and thus $f(p)\geq f(q)>f(p')$. As $f(p)$
is minimum then $p'\in K_{n-1}$ and, by Remark
\ref{Remark_symmetry_triangle_n}, $m|z_p$ as $m|w(z)_q$. The
induction step runs in a similar way to the earlier case.
\end{proof}

%% file: webb.tex
\section{Webb's conjecture}\label{section_webb}
In [7], Brown introduces the so called Brown's complex of a finite
group. Given a finite group $G$ and a prime $p$, its associated
Brown's complex $\Ss_p(G)$ is the geometrical realization of the
poset of non trivial $p$-subgroups of $G$. Webb conjectured that the
orbit space $\Ss_p(G)/G$ (as topological space) is contractible.
This conjecture was first proven by Symonds in \cite{symonds},
generalized for blocks by Barker \cite{barker1,barker2} and extended
to arbitrary (saturated) fusion systems by Linckelmann
\cite{markus}.

The works of Symonds and Linckelmann prove the contractibility of
the orbit space by showing that it is simply connected and acyclic,
and invoking Whitehead's Theorem. Both proofs of acyclicity work on
the subposet of normal chains (introduced by Knorr and Robinson
\cite{robinson} for groups). Symonds uses that the subposet of
normal chains is $G$-equivalent to Brown's complex, as was proved by
Th\'{e}venaz and Webb in \cite{G-homot}. Linckelmann proves on his
own that, also for fusion systems, the orbit space and the orbit
space on the normal chains has the same cohomology \cite[Theorem
4.7]{markus}. In this chapter we shall apply the results of Section
\ref{section_global} to prove, in an alternative way, that the orbit
space on the normal chains is acyclic.

Let $(S,\Ff)$ be a saturated fusion system where $S$ is a $p$-group
\cite{blo2}. Consider its subdivision category $\Ss(\Ff)$ (see
\cite[Section 2]{markus}) and the poset $[\Ss(\Ff)]$. An object in
$[\Ss(\Ff)]$ is an $\Ff$-isomorphism class of chains
$$
[Q_0<Q_1<\ldots<Q_n]
$$
where the $Q_i$'s are subgroups of $S$. The subcategory
$([\Ss(\Ff)]\downarrow [Q_0<\ldots<Q_n])$ has objects
$[Q_{i_0}<\ldots<Q_{i_m}]$ with $0\leq m\leq n$ and $0\leq
i_1<i_2<\ldots<i_m\leq n$ (see \cite[Section 2]{markus} again). For
example, $([\Ss(\Ff)]\downarrow [Q_0<Q_1<Q_2])$ is
$$
\xymatrix{ [Q_0]\ar[r]\ar[rd] & [Q_0<Q_1]\ar[rd]\\
[Q_1]\ar[ru]\ar[rd]&[Q_0<Q_2]\ar[r]&[Q_0<Q_1<Q_2].\\
[Q_2]\ar[r]\ar[ru]&[Q_1<Q_2]\ar[ru] }
$$

Then it is clear that $[\Ss(\Ff)]$ is a simplex-like poset.
Following Linckelmann's notation, denote by $\Ss_\lhd(\Ff)$ the full
subcategory of $\Ss(\Ff)$ which objects $Q_0<\ldots<Q_n$ with
$Q_i\lhd Q_n$ for $i=0,\ldots,n$. Also, denote by  $[\Ss_\lhd(\Ff)]$
the subdivision category of $\Ss_\lhd(\Ff)$, which is a sub-poset of
$[\Ss(\Ff)]$. By \cite[Theorem 4.7]{markus},
$H^*([\Ss(\Ff)];\Z)=H^*([\Ss_\lhd(\Ff)];\Z)$. Our goal in this
section is to prove that $H^n([\Ss_\lhd(\Ff)];\Z)=0$ for $n\geq 1$
and $H^0([\Ss_\lhd(\Ff)];\Z)=\Z$. It is straightforward that
$[\Ss_\lhd(\Ff)]$ is a simplex-like poset and thus, by Lemma
\ref{lem_coveringfamily_for_locallyDelta}, there exists an adequate
covering family for the graded poset $[\Ss_\lhd(\Ff)]^{op}$. We
shall build a global covering family $\K$ for
$[\Ss_\lhd(\Ff)]^{op}$. This family $\K$ will verify $b^\K_n=0$ for
$n\geq 1$ and $b^\K_0=1$, and so Theorem
\ref{proposition_cohomology_cochain_global_covering} will give the
desired result.

The definition of the global covering family is as follows, and it
is related with the pairing defined by Linckelmann in
\cite[Definition 5.5]{markus}. The notion of paired chains was used
by Knörr and Robinson in several forms throughout \cite{robinson}.

\begin{Defi}\label{defi_global_covering_family for normal chains}
For the graded poset $[\Ss_\lhd(\Ff)]^{op}$ define the subsets
$\K=\{K_n\}_{n\geq 0}$ by
$$
K_n={\big\{}[Q_0<\ldots<Q_n]\textrm{ $|$
}[Q_0<\ldots<Q_n]=[Q'_0<\ldots<Q'_n]\Rightarrow \cap_{i=0}^n
N_S(Q'_i)=Q'_n{\big\}}.
$$
\end{Defi}

The fact that $\K$ is defined through a pairing provides (see below)
a bijection $\psi:\Ob_n([\Ss_\lhd(\Ff)]^{op})\setminus
K_n\rightarrow K_{n+1}$ for each $n\geq 0$. This proves that
$b^\K_n=0$ for each $n\geq 1$. Moreover, $K_0=\{[S]\}$ and thus
$b^\K_0=1$. Next we prove that the family $\K=\{K_n\}_{n\geq 0}$
defined in \ref{defi_global_covering_family for normal chains} is a
global covering family for $[\Ss_\lhd(\Ff)]^{op}$. We use
terminology and results from \cite[Appendix]{blo2}.

For any chain $Q_0<\ldots<Q_n$ in $\Ss_\lhd(\Ff)$ define the
following subgroup of automorphisms of $Q_n$
$$
A_{Q_0<\ldots<Q_n}=\{\alpha\in \Aut(Q_n)|\alpha(Q_i)=Q_i,
i=0,\ldots,n\}.
$$

Then,
$$
N^{A_{Q_0<\ldots<Q_n}}_S(Q_n)=\cap_{i=0}^n N_S(Q_i).
$$
If $[Q_0<\ldots<Q_n]=[Q'_0<\ldots<Q'_n]$ then there is $\varphi\in
\Iso_\Ff(Q_n,Q'_n)$ with $Q'_i=\varphi(Q_i)$ for $i=0,\ldots,n$ and
$$
\varphi A_{Q_0<\ldots<Q_n} \varphi^{-1} = A_{Q'_0<\ldots<Q'_n}.
$$
By \cite[A.2(a)]{blo2}, $Q_n$ is fully
$A_{Q_0<\ldots<Q_n}$-normalized if and only if
$|N^{A_{Q_0<\ldots<Q_n}}_S(Q_n)|$ is maximum among
$|N^{A_{Q'_0<\ldots<Q'_n}}_S(Q_n)|$ with
$[Q'_0<\ldots<Q'_n]=[Q_0<\ldots<Q_n]$, i.e., if and only if
$|\cap_{i=0}^n N_S(Q_i)|$ is maximum among $|\cap_{i=0}^n
N_S(Q'_i)|$ with $[Q'_0<\ldots<Q'_n]=[Q_0<\ldots<Q_n]$. Notice that
in the isomorphism class of chains $[Q_0<\ldots<Q_n]$ there is
always a representantive $Q'_0<\ldots<Q'_n$ which is fully
$A_{Q'_0<\ldots<Q'_n}$-normalized, and that any two representantives
$Q'_0<\ldots<Q'_n$ and $Q''_0<\ldots<Q''_n$ of $[Q_0<\ldots<Q_n]$
which are fully $A_{Q'_0<\ldots<Q'_n}$-normalized and fully
$A_{Q''_0<\ldots<Q''_n}$-normalized respectively verify
$$
|\cap_{i=0}^n N_S(Q'_i)|=|\cap_{i=0}^n N_S(Q''_i)|.
$$

Thus, Definition \ref{defi_global_covering_family for normal chains}
is equivalent to
\begin{Defi}
For the graded poset $[\Ss_\lhd(\Ff)]^{op}$ define the subsets
$\K=\{K_n\}_{n\geq 0}$ by
$$
K_n=\{[Q'_0<\ldots<Q'_n]\text{ $|$ $Q'_n$ fully
$A_{Q'_0<\ldots<Q'_n}$-normalized and $\cap_{i=0}^{n}
N_S(Q'_i)=Q'_n$}\}.
$$
\end{Defi}

\begin{Lem}\label{lemma_bijection_K_n}
For any $n\geq 0$ there is a bijection
$$
\psi:\Ob_n([\Ss_\lhd(\Ff)]^{op})\setminus K_n\rightarrow K_{n+1}.
$$
\end{Lem}
\begin{proof}
Take $[Q_0<\ldots<Q_n]\in \Ob_n([\Ss_\lhd(\Ff)]^{op})\setminus K_n$
and a representantive $Q'_0<\ldots<Q'_n$ which is fully
$A_{Q'_0<\ldots<Q'_n}$-normalized. Then $\cap_{i=0}^{n}
N_S(Q'_i)>Q'_n$. Define
$\psi([Q_0<\ldots<Q_n])=[Q'_0<\ldots<Q'_n<\cap_{i=0}^{n}
N_S(Q'_i)]$. The proof is divided in four steps:

\textbf{a) $\psi$ is well defined.} Take another representantive
$Q''_0<\ldots<Q''_n$ which is fully
$A_{Q''_0<\ldots<Q''_n}$-normalized. Then, by \cite[A.2(c)]{blo2},
there is a morphism
$$
\varphi\in
\Hom_\Ff(N_S^{A_{Q'_0<\ldots<Q'_n}}(Q'_n),N^{A_{Q''_0<\ldots<Q''_n}}_S(Q'_n))
$$
with $\varphi(Q'_i)=Q''_i$ for $i=0,\ldots,n$. As $|\cap_{i=0}^{n}
N_S(Q'_i)|=|\cap_{i=0}^{n} N_S(Q''_i)|$ then $\varphi$ is an
isomorphism onto $N^{A_{Q''_0<\ldots<Q''_n}}_S(Q'_n)$ and thus
$$
[Q'_0<\ldots<Q'_n<\cap_{i=0}^{n}
N_S(Q'_i)]=[Q''_0<\ldots<Q''_n<\cap_{i=0}^{n} N_S(Q''_i)].
$$

\textbf{b) $\psi([Q_0<\ldots<Q_n])$ belongs to $K_{n+1}$.} We have
$\psi([Q_0<\ldots<Q_n])=[Q'_0<\ldots<Q'_n<\cap_{i=0}^{n} N_S(Q'_i)]$
where $Q'_0<\ldots<Q'_n$ is fully $A_{Q'_0<\ldots<Q'_n}$-normalized.
Take any representantive $Q''_0<\ldots<Q''_n<Q''_{n+1}$ in
$[Q'_0<\ldots<Q'_n<\cap_{i=0}^{n} N_S(Q'_i)]$. If it were the case
that $Q''_{n+1}< \cap_{i=0}^{n+1} N_S(Q''_i)$ then we would have
$$
\cap_{i=0}^{n} N_S(Q'_i)\cong Q''_{n+1} < \cap_{i=0}^{n+1}
N_S(Q''_i)\leq  \cap_{i=0}^{n} N_S(Q''_i),
$$
which is in contradiction with $Q'_0<\ldots<Q'_n$ being fully
$A_{Q'_0<\ldots<Q'_n}$-normalized.

\textbf{c) $\psi$ is injective.} Suppose we have $[Q_0<\ldots<Q_n]$
and $[R_0<\ldots<R_n]$ with
$$
[Q'_0<\ldots<Q'_n<\cap_{i=0}^{n}
N_S(Q'_i)]=[R'_0<\ldots<R'_n<\cap_{i=0}^{n} N_S(R'_i)].
$$
Then
$[R_0<\ldots<R_n]=[R'_0<\ldots<R'_n]=[Q'_0<\ldots<Q'_n]=[Q_0<\ldots<Q_n]$.

\textbf{d) $\psi$ is surjective.} Take $[Q_0<\ldots<Q_n<Q_{n+1}]$ in
$K_{n+1}$. We check that
$$\psi([Q_0<\ldots<Q_n])=[Q_0<\ldots<Q_n<Q_{n+1}].$$
Take a representantive $Q'_0<\ldots<Q'_n$ in $[Q_0<\ldots<Q_n]$
which is fully $A_{Q'_0<\ldots<Q'_n}$-normalized. Then
$[Q_0<\ldots<Q_n]\in \Ob_n([\Ss_\lhd(\Ff)]^{op})\setminus K_n$ and
$\psi([Q_0<\ldots<Q_n])=[Q'_0<\ldots<Q'_n<\cap_{i=0}^{n}
N_S(Q'_i)]$. Then, by \cite[A.2(c)]{blo2}, there is $$\varphi\in
\Hom_\Ff(N_S^{A_{Q_0<\ldots<Q_n}}(Q_n),N^{A_{Q'_0<\ldots<Q'_n}}_S(Q_n))
$$
with $\varphi(Q_i)=Q'_i$ for $i=0,\ldots,n$. As
$Q_{n+1}=\cap_{i=0}^{i=n+1} N_S(Q_i)$ then $Q_{n+1}\leq
\cap_{i=0}^{n} N_S(Q_i)$ and $\varphi(Q_{n+1})\leq \cap_{i=0}^{n}
N_S(Q'_i)$. If it were the case that $\varphi(Q_{n+1})<
\cap_{i=0}^{n} N_S(Q'_i)$ then we would have
$$
\varphi(Q_{n+1})< N_{\cap_{i=0}^{n}
N_S(Q'_i)}(\varphi(Q_{n+1}))=\cap_{i=0}^{n+1}
N_S(\varphi(Q_i))=\varphi(\cap_{i=0}^{n+1}
N_S(Q_i))=\varphi(Q_{n+1}),
$$
a contradiction. Thus, $\varphi(Q_{n+1})=\cap_{i=0}^{n} N_S(Q'_i)$
and the proof is finished.
\end{proof}

\begin{Lem}
The family $\K=\{K_n\}_{n\geq 0}$ defined in
$\ref{defi_global_covering_family for normal chains}$ is a global
covering family for $[\Ss_\lhd(\Ff)]^{op}$.
\end{Lem}
\begin{proof}
We start proving that for any $n\geq 0$ the map
$$
\liminv F_n\rightarrow \prod_{i\in K_n} F_n(i)
$$
is a monomorphism. Take $\psi\in \liminv F_n$ such that $\psi(i)=0$
for each $i\in K_n$. If there is no object of degree greater than
$n$ then $K_n=\Ob_n([\Ss_\lhd(\Ff)]^{op})$ and we are done. If not,
we prove that $\psi(j)=0$ for each $j$ of degree $n+1$ by induction
on $|Q_{n+1}|$. This is enough to see that $\psi$ is zero as $F_n$
is $n$-condensed. The base case is $j=[Q_0<\ldots<Q_{n+1}]$ with
$|Q_{n+1}|$ maximal. This implies that $J^j_n\subseteq K_n$. Then
$\psi(j)$ goes to zero by the monomorphism
$$
F_n(j)\rightarrow \prod_{i\in J^{j}_n} F_n(i),
$$
and thus $\psi(j)=0$. For the induction step consider
$j=[Q_0<\ldots<Q_{n+1}]$ and
$j'=[Q_0<\ldots<\widehat{Q_l}<\ldots<Q_{n+1}] \in J^j_n$ with $0\leq
l<n$. Then, either $j'\in K_n$ and $\psi(j')=0$, or $j'\notin K_n$
and there is an arrow in $[\Ss_\lhd(\Ff)]^{op}$
$$
j''=[Q'_0<\ldots<\widehat{Q'_l}<\ldots<Q'_{n+1}< \cap_{i=0,i\neq
l}^{m} N_S(Q'_i)]\rightarrow
j'=[Q_0<\ldots<\widehat{Q_l}<\ldots<Q_{n+1}].
$$
In the latter case $\psi(j'')=0$ by the induction hypothesis, and
thus $\psi(j')=0$ too. As before, since the map $F_n(j)\rightarrow
\prod_{i\in J^{j}_n} F_n(i)$ is a monomorphism, $\psi(j)=0$.

Now we prove that for any $n\geq 1$ the map
$$
\omega: \prod_{i\in \Ob_{n-1}(\P)\setminus K_{n-1}}
F_{n-1}(i)\rightarrow \prod_{i\in K_n} F_n(i)
$$
is pure. Take $y\in \prod_{i\in K_n} F_n(i)$, $m\geq 1$ and $x\in
\prod_{i\in \Ob_{n-1}(\P)\setminus K_{n-1}} F_{n-1}(i)$ with
$$
m\cdot y = \omega(x).
$$
We want to find $x'$ with $m\cdot x'=x$. We prove that $x_i$ is
divisible by $m$ for each $i=[Q_0<\ldots<Q_{n-1}]\in
\Ob_{n-1}(\P)\setminus K_{n-1}$ by induction on $|Q_{n-1}|$. The
base case is $i=[Q_0<\ldots<Q_{n-1}]$ with $|Q_{n-1}|$ maximal.
Consider the arrow in $[\Ss_\lhd(\Ff)]^{op}$
$$
j=[Q'_0<\ldots<Q'_{n-1}< \cap_{i=0}^{n-1} N_S(Q'_i)]\rightarrow
[Q_0<\ldots<Q_{n-1}].
$$
As $Q_{n-1}$ is maximal then $J^j_{n-1}\subseteq K_{n-1}$. Then
$m\cdot y_j=\omega(x)_j$ is the image of $(x_i,0,\ldots,0)$ by the
map
$$
\ker_{F_{n-1}}'(j)=\prod_{l\in (j\downarrow
{[\Ss_{\lhd}(\Ff)]})_{n-1}} F_{n-1}(i)\stackrel{\pi_j}\rightarrow
F_n(j).
$$
As $F_n(j)\cong \prod_{l\in (j\downarrow {[\Ss_{\lhd}(\Ff)]}
)_{n-1}\setminus J^j_{n-1}} F_{n-1}(l)=F_{n-1}(i)$ by Remark
\ref{rmk_section} then $m$ divides $x_i$.

For the induction step consider $i=[Q_0<\ldots<Q_{n-1}]\in
\Ob_{n-1}(\P)\setminus K_{n-1}$ and $j=[Q'_0<\ldots<Q'_{n-1}<
\cap_{i=0}^{i=n-1} N_S(Q'_i)]$. As before, $m\cdot y_j=\omega(x)_j$
is the image of $\tilde{x}=x|{(j\downarrow
{[\Ss_{\lhd}(\Ff)]})_{n-1}}$ by the map
$$
\ker_{F_{n-1}}'(j)=\prod_{l\in (j\downarrow
{[\Ss_{\lhd}(\Ff)]})_{n-1}} F_{n-1}(i)\stackrel{\pi_j}\rightarrow
F_n(j).
$$
By Remark \ref{rmk_section} again,
$$
m\cdot y_j=\pi_j(\tilde{x}-(\lambda_j\circ
s_j)(\tilde{x}))=\pi_j((x_i-(\lambda_j\circ
s_j)(x|{J^j_n}),0,\ldots,0)).
$$
Now, by the induction hypothesis, for each $l\in J^j_n$ either $l\in
K_{n-1}$ and $x_l=0$, either $l\notin K_{n-1}$ and thus, by the
induction hypothesis, $m$ divides $x_l$. Then $m$ divides
$x|{J^j_n}$ and so, by the equation above and the isomorphism
$F_n(j)\cong F_{n-1}(i)$, $m$ divides $x_i$ too.
\end{proof}

%% file: coxeter.tex
\section{Coxeter groups}\label{section_coxeter}
It is well known that the Coxeter complex associated to a finite
Coxeter group $(W,S)$ is a simplicial decomposition of a sphere of
dimension $|S|-1$, and thus has the homotopy type of a sphere. In
case $W$ is infinite then this complex becomes contractible. In this
section we prove that the cohomology of the Coxeter complex is that
of a sphere if $W$ is finite or trivial if $W$ is infinite. We do it
by using the techniques of earlier sections. More precisely, we
construct a global covering family $\K$ for the Coxeter complex with
$b^K_0=b^K_{|S|-1}=1$ and $b^K_n=0$, $n\neq 0, |S|-1$ for $W$
finite, and with $b^K_0=1$ and $b^K_n=0$, $n\neq 0$ for $W$
infinite. We will use general facts about Coxeter groups which can
be found in \cite{bourbaki}, \cite{humphreys} or \cite{garrett}, as
well as we will recall some basic definitions and statements.

Let $(W,S)$ be a Coxeter system where $W$ is a Coxeter group and
$S=\{s_1,...,s_N\}$ is a set of generators (which we always assume
is finite). For any word $w\in W$ its length $\l(w)$ is the minimum
number of generators from $S$ that are needed to write it. If $W$ is
finite there is a unique element of maximal length which we denote
by $w_0$. For every subset $I\subset S$ we have the parabolic
subgroup $W_I\leq W$ generated by the generators belonging to $I$.
We also define the following subset
$$
W^I=\{w\in W| \l(ws)>\l(w)\text{ for all $s\in I$}\}.
$$
In \cite[5.12]{humphreys} and \cite[excersise 3, p. 37]{bourbaki} is
proven the following:
\begin{Lem}
Fix $I\subseteq S$. Given $w\in W$, there is a unique $u\in W^I$ and
a unique $v\in W_I$ such that $w=uv$. Then $\l(w)=\l(u)+\l(v)$.
Moreover, $u$ is the unique element of smallest length in the coset
$wW_I$.
\end{Lem}

From this lemma it is straightforward that there is a bijection
\begin{equation}\label{equ_coxeter_bij}
W/W_I\rightarrow W^I
\end{equation}
which sends the coset $wW_I\in W/W_I$ to the unique element $u\in
W^I$ of smallest length in the coset $wW_I$. Notice that $uW_I=wW_I$
for $w$ and $u$ as in the lemma. Next we describe the Coxeter
complex associated to $(W,S)$: it is a simplicial complex which
simplices are the cosets of proper parabolic subgroups with
inclusion reversed:
$$
wW_I\leftarrow w'W_{I'}\Leftrightarrow I\subset I'\text{ and
}w^{-1}w'\in W_{I'}.
$$
The dimension of the simplex $wW_I$ is $n=|S|-|I|-1$. Thus, vertices
correspond to maximal subsets $I$ of $S$ and facets, i.e., maximal
dimensional faces, correspond to $I=\emptyset$. Notice that it is a
pure simplicial complex. By the bijection (\ref{equ_coxeter_bij}) we
can write the coset $wW_I$ as $uW_I$ with $u\in W^I$, where $w$ and
$u$ are as in the previous lemma. We will do this in the rest of the
section. For any $u\in W$ define (cf. \cite[1.11]{humphreys})
$$
S_u=\{s\in S|\l(us)>\l(u)\}
$$
and let $s_u=\max\{S_u\}$ be the maximum element in $S_u$ with
respect to the order in $S$: $s_1<...<s_N$. Notice that $S_u=S$ if
and only if $u=1$. Moreover, if $W$ is finite then $S_u=\emptyset$
if and only $u=w_0$, the unique word of maximal length. It is also
clear that
\begin{equation}\label{equ_Hump_1.11_heart}
\text{$u\in W^I\Leftrightarrow I\subseteq S_u$.}
\end{equation}

Before defining the global covering family notice that, by Lemma
\ref{lem_coveringfamily_for_locallyDelta}, there is an adequate
local covering family for $\P^{op}$. Now, for each $n\geq 0$ define
$$
K_n=\{uW_I|\text{$ |I|=|S|-1-n$ and $s_u\notin I$}\}.
$$
The condition $|I|=|S|-1-n$ in the definition above just states that
$uW_I$ correspond to a simplex of dimension $n$.
\begin{Lem}
The family $\K=\{K_n\}_{n\geq 0}$ is a global covering family for
$\P^{op}$. Moreover, $b^\K_0=1$. If $W$ is finite then
$b^\K_{|S|-1}=1$ and $b^\K_n=0$ for $n\neq 0,|S|-1$. If $W$ is
infinite then $b^\K_n=0$ for $n\neq 0$.
\end{Lem}
By Theorem \ref{proposition_cohomology_cochain_global_covering} this
lemma imply that the integral cohomology of $\P$ is that of a sphere
of dimension $|S|-1$ if $W$ is finite, and that of a point if $W$ is
infinite.
\begin{proof}
First we compute the numbers $b^K_n$ and afterwards we will prove
$\K$ satisfies Definition \ref{defi_global covering family}. Fix
$u\in W$ and consider the contribution which it makes to $\P^{op}$.
By Equation (\ref{equ_Hump_1.11_heart}) $u\in W^I$ if and only if
$I\subseteq S_u$. Thus, denoting by $\P_u$ the sub-poset with
objects $\{uW_I, I\subseteq S_u\}$ we have that $\P$ is contained in
the disjoint union $\bigcup_{u\in W} \P_u$. It is a combinatorial
exercise (cf. Lemma \ref{lem_coveringfamily simplex}) that there is
a bijection between $K_{n+1}\cap \P_u$ and $\{\P_n\setminus
K_n\}\cap \P_u$ for $n\geq 0$ if $u\neq w_0$. This gives, in case
$W$ is infinite, $b^K_n=0$ for $n\geq 1$. If $W$ is finite then
$\P_{w_0}=\{w_0W_\emptyset\}$ and $w_0W_\emptyset$ is a facet which
belongs to $K_{|S|-1}$. Then we obtain $b^K_{|S|-1}=1$ and $b^K_n=0$
for $|S|-1> n\geq 1$. Finally, it is an easy consequence of the
definition that $K_0=\{1\cdot W_{S\setminus\{s_N\}}\}$ and thus
$b^K_0=1$.

Now fix $n\geq 0$. We show that the restriction map
$$
\liminv F_n\rightarrow \prod_{uW_I\in K_n} F_n(uW_I)
$$
is a monomorphism, where $F_n:\P^{op}\rightarrow \Ab$ are the
functors obtained from Lemma \ref{cohomology with adequate local
covering family} applied to $\P^{op}$.

Take $\psi\in \liminv F_n=\hom_{\Ab^{\P}}(c_{\Z},F_n)$ which is
mapped to zero by the restriction map. To prove that $\psi=0$ it is
enough to prove that $\psi(uW_I)=0$ for each simplex $uW_I$ of
dimension $n+1$ (as in the proof of Theorem
\ref{morse_prop_modelo}). We do this by induction on the length
$\l(u)$.

We start with $uW_I$ with $l(u)=0$, i.e., $u=1$. Then $S_u=S$. If
$s_u=s_N\notin I$, i.e., $u\in K_{n+1}$, then there $n+1$ $n$-faces
of $u$ which are in $K_n$. These are the cosets $W_J$ where $J=I\cup
\{s\}$ with $s\in S\setminus\{I\cup\{s_N\}\}$. The remaining face,
i.e., $W_J$ with $J=I\cup\{s_N\}$, is not in $K_n$. By Remark
\ref{Remark_symmetry_triangle_n}, $\psi(W_I)=0$. Now assume that
$s_u=s_N\in I$, i.e., $u\notin K_{n+1}$. Consider any of the $n+2$
$n$-faces $W_J$ of $W_I$, where $J=I\cup \{s\}$ with $s\in
S\setminus I$. Then $W_J$ is also a face of then $n+1$ simplex
$W_{I\setminus\{s_N\}\cup\{s\}}$. By the preceding argument for the
case $W_{I\setminus\{s_N\}\cup\{s\}}\in K_{n+1}$ we obtain
$\psi(W_{I\setminus\{s_N\}\cup\{s\}})=0$, and thus $\psi(W_J)=0$.
Then $\psi$ is zero in all the faces of $W_I$. By Remark
\ref{Remark_symmetry_triangle_n} again we obtain $\psi(W_I)=0$.

Next we do the induction step: take an $n+1$ dimensional simplex
$uW_I$ with $l(u)>0$. First assume that $u\in K_{n+1}$, i.e.,
$s_u\notin I$. The faces of $uW_I$ are the $n+2$ $n$-simplices
$uW_J$ where $J=I\cup\{s\}$ and $s\in S\setminus I$. Notice that, as
$I\subseteq S_u$,  $S\setminus I=S_u\setminus I \cup S\setminus
S_u$, where the union is disjoint. Take first $s\in S\setminus S_u$.
Then $l(us)<l(u)$ and, as $us\in uW_J$, the unique element $u'$ of
minimal length in $uW_J$ is different from $u$ and has smaller
length. Then $uW_J=u'W_J$ is also a face of $u'W_I$. By the
induction hypothesis and because $l(u')<l(u)$ we have
$\psi(u'W_I)=0$ and thus $\psi(uW_J)=0$. Now take $s\in S_u\setminus
I$. Then $uW_J\in K_n$ unless $s=s_u$. Then $\psi$ is zero in all
but one face of $uW_I$ and, by Remark
\ref{Remark_symmetry_triangle_n}, $\psi(uW_I)=0$.

Now assume that $u\notin K_{n+1}$, i.e., $s_u\in I$, and take a face
$uW_J$ as before. If $J=I\cup\{s\}$ with $s\in S\setminus S_u$ then
arguing as before we obtain that $\psi(uW_J)=0$. If $J=I\cup\{s\}$
with $s\in S_u\setminus I$ then $uW_J$ is also a face of the $n+1$
simplex $uW_{I\setminus\{s_N\}\cup\{s\}}$. By the preceding argument
for the case $uW_{I\setminus\{s_N\}\cup\{s\}}\in K_{n+1}$ we have
$\psi(uW_{I\setminus\{s_N\}\cup\{s\}})=0$, and thus $\psi(W_J)=0$.
Then $\psi$ takes the value zero in all the faces of $uW_I$ and by
Remark \ref{Remark_symmetry_triangle_n} $\psi(uW_I)=0$.

Finally, fix $n\geq 1$. The proof that the map
$$
w:\prod_{uW_I\in \P_{n-1}\setminus K_{n-1}} F_{n-1}(uW_I)\rightarrow
\prod_{u'W_{I'}\in K_n} F_n(u'W_{I'})
$$
is pure is made as in Theorem \ref{morse_prop_modelo}. It uses
induction on the length $l(u)$ for $(n-1)$-dimensional simplices
$uW_I$, and the base case is $u=1$.
\end{proof}